\documentclass[a4paper,11pt,reqno]{amsart}

\usepackage{listings}

\usepackage{graphicx}
\usepackage{amsmath}
\usepackage{amssymb}
\usepackage{amsfonts}
\usepackage[hang]{caption}
\usepackage{comment}
\usepackage{mathrsfs}
\usepackage{graphics}
\usepackage{float}
\usepackage{mathtools}

\usepackage{rotating}

\usepackage[utf8]{inputenc}

\DeclareFixedFont{\ttb}{T1}{txtt}{bx}{n}{8} 
\DeclareFixedFont{\ttm}{T1}{txtt}{m}{n}{8}  

\usepackage{color}
\definecolor{deepblue}{rgb}{0,0,0.85}
\definecolor{deepred}{rgb}{0.6,0,0}
\definecolor{deepgreen}{rgb}{0,0.5,0}
\definecolor{deepbrown}{rgb}{0.88,0.52,0.26}
\definecolor{lavender}{rgb}{0.81,0.41,0.83}

\usepackage{listings}

\newcommand\pythonstyle{\lstset{
language=Python,
basicstyle=\ttm,
morekeywords={str},                
keywordstyle=\ttb\color{deepblue},
emph={MyClass,__init__},          
emphstyle=\ttb\color{deepred},    
stringstyle=\color{deepgreen},
commentstyle = \color{deepgreen},
frame=tb,                         
showstringspaces=false,
classoffset=2,
morekeywords={self},
keywordstyle=\color{deepbrown},
}}

\lstnewenvironment{python}[1][]
{
\pythonstyle
\lstset{#1}
}
{}

\newcommand\pythoninline[1]{{\pythonstyle\lstinline!#1!}}

\graphicspath{{eps/}}

\newcommand{\RomanNumeralCaps}[1]
    {\MakeUppercase{\romannumeral #1}}

\makeatletter

\newtheorem{theorem}{Theorem}[section]

\@addtoreset{equation}{section}

\newtheorem{remark}[theorem]{Remark}
\newtheorem{proposition}[theorem]{Proposition}

\makeindex

\begin{document}

\title[Extending Snow's algorithm]
{Extending Snow's algorithm for computations in the finite Weyl groups}
         \author{Rafael Stekolshchik}

\date{}

\begin{abstract}
   In 1990, D.Snow proposed an effective algorithm for computing
   the orbits of finite Weyl groups.
   Snow's algorithm is designed for computation of weights, $W$-orbits and elements of the Weyl group.
   An extension of Snow's algorithm is proposed, which allows 
   to find pairs of mutually inverse elements together with the calculation of $W$-orbits
   in the same runtime cycle.
   This  simplifies the calculation of conjugacy classes in the Weyl group.
   As an example, the complete list of elements of the Weyl group $W(D_4)$
   obtained using the extended Snow's algorithm is given.
   The elements of $W(D_4)$ are specified in two ways:
   as reduced expressions and as matrices of the faithful representation.
   We present a partition of this group into conjugacy classes with elements specified
   as reduced expressions.
   Various forms are given for representatives of the conjugacy classes of $W(D_4)$:
   using Carter diagrams, using reduced expressions and using signed cycle-types.
   In the appendix, we provide an implementation of the algorithm in Python.
\end{abstract}

\maketitle

\section{\bf Introduction}

\subsection{Snow's algorithm: finding $W$-orbits}
In 1990, D.Snow in \cite{Sn90} proposed an effective algorithm for computing
the orbits of the finite Weyl groups. The algorithm
starts with a certain dominant weight and acts on it
by all the simple reflections. This operation produces the complete list of weights of level $1$ and
the complete list of all elements of length $1$ in the Weyl group $W$.
At the next step, we use reflections again to get
a list of weights of level 2 and all elements  of length $2$, and so on.
This approach has a repetition problem: the same weight can be obtained in several ways,
and the list of elements of the Weyl group lying in some level will contain duplicate elements.
Snow presented a solution showing which weight $v$ should be taken at the level $L_k$
and which reflection $s_i$ should be applied to $v$ to get the given weight $\xi$
at the level $L_{k+1}$. By Snow's algorithm, choosing  $v$ and $s_i$ can be carried out by the unique way.
This solution avoids duplicate elements, see Section \ref{sec_avoid_rep}.

The computation of the elements of the Weyl group in Snow's algorithm is based on the following fact:
there is a one-to-one correspondence between the Weyl chambers and the elements of the Weyl group,
and the Weyl group acts transitively on the set of Weyl chambers.
Each element from the closure of the fundamental Weyl chamber generates
a Weyl group orbit ($W$-orbit) whose length coincides with
the order of the Weyl group.  The $W$-orbit is constructed
under the action of the Weyl group on some dominant weight.
The  weights of the $W$-orbit are constructed together with the elements of the
Weyl group $W$.

Let $\varPhi$ be the root system associated with a certain semisimple Lie algebra $\mathcal{L}$,
and $\mathcal{E}$ be a real space spanned by the roots of $\varPhi$.
A {\it weight} is an element $\xi \in \mathcal{E}$ such that $\langle \xi, \alpha \rangle \in \mathbb{Z}$
for all roots $\alpha \in \varPhi$.
The set of weights $\Lambda$ forms a subgroup of $\mathcal{E}$, i.e.,
$\varPhi \subset \Lambda \subset \mathcal{E}$.
The  significance of the  weights theory is largely determined by the highest weight theorem
in the representation theory of semisimple Lie algebras\footnotemark[1].
\footnotetext[1]{The highest weight theorem was proved by E.Cartan in $1913$, \cite{C13},
see $\S$\ref{sec_highest}.}

Snow's algorithm produces the weights of $W$-orbits and elements of the Weyl group by {\it levels}.
For any $\xi \in \mathcal{E}$ there exist $w \in W$ and $v$
from the closure $\overline{C}$ of the fundamental Weyl chamber such that $\xi = w(v)$,
see Theorem \ref{fact_1}.  The level of $\xi$ is as follows:
\begin{equation}
   \text{level}(\xi) = l(w),
\end{equation}
where $l(w)$ is the length of $w$ given as a reduced expression.
The {\it level of weight} $\xi$  is equal to the number of reflections needed to move $\xi$
to some dominant weight lying in the closure of the fundamental chamber $\overline{C}$,
see Proposition \ref{prop_inique_dominant}.

In Section \ref{sec_Snow}, we will look at some of the details of Snow's algorithm.
The sizes of all levels and the total computation time for cases
$B_7$, $D_8$, $E_7$, $B_8$ are gathered in Table \ref{tab_levels_BDE}.

\subsection{Extended Snow's algorithm: finding inverse elements}
To construct conjugacy classes of a group,
one must first find all pairs of mutually inverse elements of the group.
In the case of the Weyl group, each element and its inverse belong to the same level.
However, even searching within a level can be quite an expensive task,
especially for very large levels, see Table \ref{tab_levels_BDE},
where the  level lengths are in the hundreds of thousands of elements.
Let
\begin{equation}
  \label{dir_order}
   w = s_{i_1}s_{i_2}\dots s_{i_{k-1}}s_{i_k}
\end{equation}
be an element of the level $L_k$.
We can find the inverse element $w^{-1}$ by reversing the order of the
reduced expression $w$:
\begin{equation}
  \label{rev_order}
   w^{-1} = s_{i_k}s_{i_{k-1}}\dots s_{i_2}s_{i_1}.
\end{equation}
However, the inverse element must be found in accordance to the repetition prevention mechanism
from Theorem \ref{th_snow}. Then the reduced expression may differ from the reverse order of $w$.

An extension of Snow's algorithm is designed to get around this obstacle:
for any element $w \in W$, one must obtain the inverse element $w^{-1}$, but this must be done
in the order specified by Theorem \ref{th_snow}. The reduced expression of
the calculated inverse element will not necessarily be of the form \eqref{rev_order}.
Bypassing the specified obstacle
achieved through the exchange of information between any element and its inverse
during the traversal performed by Snow's algorithm. This information  exchange
is carried out using the dictionary mechanism described in Section \ref{sec_snow_ext}.

The Weyl group $W(D_4)$ contains $192$ elements.
In Section \ref{sec_levels_D4}, all elements of $W(D_4)$  are divided into $12$ levels.
The elements of $W(D_4)$
are specified in two ways: as matrices and as reduced expressions,
see Tables \ref{tab_levels_012} - \ref{tab_level_11_12}.
For each element $w$, we provide
also the reduced expression of the inverse element and its location.

The partition of the group $W(D_4)$ into conjugate classes is given in Section \ref{sec_conj_classes}.
There are $13$ conjugacy classes including the trivial
class containing only identity element $e$, see Tables \ref{tableCCL_0123} - \ref{tableCCL_12}.
For each element $w$ of the conjugacy class, we provide
the level number $k$ such that $w \in L_k$ and the position of $w$ in the level $L_k$.
This information allows to find the element $w$ in the tables of levels of Section \ref{sec_levels_D4}.

The execution time of the extended Snow's algorithm
for Weyl groups $B_7$, $D_8$, $E_7$, $B_8$  on CPU 3.7 GHz/Python 3.7.3  are as follows:
\begin{equation*}
  \begin{array}{lrrr}
     & B_7  & 645120   \text{ elements }  & 59 \text{ sec } \\
     & E_7  & 2903040  \text{ elements }  & 269 \text{ sec } \\
     & D_8  & 5169960  \text{ elements }  & 570 \text{ sec } \\
     & B_8  & 10321920 \text{ elements }  & 1153 \text{ sec }
  \end{array}
\end{equation*}
For the execution time for each level, see Table \ref{tab_levels_BDE}.

Appendix \ref{app_sect_weights} lists some properties of weights related to Lie algebras and Weyl groups.
An implementation of the extended Snow's algorithm in Python is given in Appendix \ref{app_sect_implem}.

\section{\bf Snow's algorithm: computation of $W$-orbits and levels}
  \label{sec_Snow}
\subsection{Computation of the $W$-orbits}

Snow's algorithm starts with a certain dominant weight and acts on it with all simple
reflections. This produces all the weights of level $1$ and a list of all elements of length $1$
in $W$. Further, we apply this procedure again and, if we ignore duplicates,
we obtain the weights of level $2$ and a required  list of elements of length $2$ in $W$.
By repeating this procedure, we compute a list of weights of any level, and the entire group $W$
can be generated if an appropriate initial weight is chosen.

\subsection{Computation of $\text{level}(\xi)$}

The algorithm provides a simple criterion for adding an orbit element to the list of weights.
Let $\xi = (x_1, \dots, x_n)$ be any weight in the basis consisting of fundamental dominant weights, see
Section \eqref{sec_fundam_weights}. What is the level of $s_i(\xi)$ for any simple reflection $s_i$?

Let $w$ be the element in $W$ such that $\xi = w(v)$ for some
$v$ from the fundamental domain $\overline{C}$ with $\text{level}(\xi) = l(w)$.
By definition of the fundamental weights \eqref{fundam_weights}, we have
\begin{equation}
  \label{coord_xi}
   \xi = \sum\limits_i x_i \bar\omega_i,  \text{ ~and~ } \\
   x_i = \langle \xi, \alpha_i \rangle = \langle w(v), \alpha_i \rangle.
\end{equation}
By \eqref{eq_def_prod} the sign of $x_i$ coincides with the sign of $(w(v), \alpha_i)$, then
\begin{equation}
 \label{signs_of_xi}
  \begin{cases}
     & x_i = 0 \quad   \Longrightarrow \quad  s_i(\xi) = \xi,  \\
     & x_i > 0 \quad \Longrightarrow \quad (w(v), \alpha_i) > 0, \\
     & x_i < 0  \quad \Longrightarrow \quad (w(v), \alpha_i) < 0.
  \end{cases}
\end{equation}
Here, the first line \eqref{signs_of_xi} follows from \eqref{eq_action_si}.
Thus, in the case of $x_i = 0$, the reflection $s_i$ does not change the level:
\begin{equation}
 \label{case_xi_is_0}
  x_i = 0 \quad  \Longrightarrow \quad \text{level}(s_i(\xi)) = \text{level}(\xi).
\end{equation}

Further, since the Cartan-Killing form is invariant under the Weyl group $W$, we have
\begin{equation}
 \label{xi_not_0}
  \begin{array}{cc}
     & x_i > 0 \quad \Longrightarrow \quad (v, w^{-1}(\alpha_i)) > 0,  \\
     & x_i < 0 \quad \Longrightarrow \quad (v, w^{-1}(\alpha_i)) < 0.   \\
  \end{array}
\end{equation}
Since $v$ is a dominant weight, we have $\langle v, \alpha \rangle \geq 0$ for all $\alpha \in \varPhi$,
see Section \eqref{sec_dominant}.
Then by Theorem \ref{th_length_plus_1}, we have
\begin{equation}
  \begin{array}{ccccc}
     & x_i > 0 & \Longrightarrow    w^{-1}(\alpha_i) \in \varPhi^+ & \Longrightarrow
               & l(s_i w) = l(w) + 1, \\
     & x_i < 0 & \Longrightarrow    w^{-1}(\alpha_i) \in \varPhi^- & \Longrightarrow
               & l(s_i w) = l(w) - 1. \\
  \end{array}
\end{equation}
Thus the level is updated as follows:
\begin{equation}
 \label{eq_3_cases}
 \text{level}(s_i(\xi) ) =
  \begin{cases}
     & \text{level}(\xi) + 1 \quad \text{ if } x_i > 0, \\
     & \text{level}(\xi)     \;\;\; \qquad \text{ if } x_i = 0, \\
     & \text{level}(\xi) - 1 \quad  \text{ if } x_i < 0.
  \end{cases}
\end{equation}

\subsection{Arranging the weights by levels}

We start from a dominant weight $\mu \in \Lambda^+$, see eq. \eqref{eq_domin_weights}.
Let $L_k$ be the $k$th level of $W\cdotp\mu$, i.e.,
\begin{equation*}
  L_k = \{\text{weights } \xi \in W\cdotp\mu \mid \text{level}(\xi) = k \}.
\end{equation*}
Then, the orbit $W\cdotp\mu$ is the disjoint union of all levels:
\begin{equation*}
   W\cdotp\mu = \bigsqcup\limits_{i = 0}^N L_i,
\end{equation*}
where $N$ is the maximal possible level in $W\cdot\mu$.
By Proposition \ref{prop_elem_max_len} the number $N$ is the number of positive roots
in $C$, since this is the maximal length of a Weyl group element.

To construct level $L_{k+1}$ from the previously computed level $L_k$, we apply reflections $s_i$.
By \eqref{eq_3_cases}, if $x_i > 0$ only reflection $s_i$ move $\xi$ from $L_k$  to $L_{k+1}$:
\begin{equation}
  \label{eq_next_levl}
  L_{k+1} = \{ s_i(\xi) ~\mid~ i = 1,\dots, l,
  ~\xi = (x_1,\dots, x_l) \in L_k, ~x_i > 0 \}.
\end{equation}

\subsection{Snow's solution to the repetition problem}
  \label{sec_avoid_rep}

\subsubsection{An example of the repetition problem}
  \label{ex_repet_problem}
  We start with the dominant weight $\lambda_0=[1,1,1,1]$ and act on this weight
  by two different elements from level 2 of the Weyl group $W(D_4)$:
  $w_1 = s_2s_3$ (Table \ref{tab_levels_012}, element 6) and $w_3 = s_3s_2$
  (Table \ref{tab_levels_012}, element 4).
  The images $\lambda_1 = w_1(\lambda_0)$ and  $\lambda_2 = w_2(\lambda_0)$ are as follows:
\begin{equation*}
  \lambda_1 = [3, -2, 1, 3] , \quad \lambda_2 = [2, 1, -2, 2].
\end{equation*}

We act by reflection $s_3$ on the weight $\lambda_1$. Here $m_3 = 1$, $\overline{c}_3 = [0, -1, 2, 0]$.
Similarly, we act by reflection $s_2$ on $\lambda_2$, where  $m_2 = 1$, $\overline{c}_2 = [-1, 2, -1, -1]$,
see formula \eqref{eq_action_si}.
\begin{equation*}
  \begin{split}
   & s_3(\lambda_1) = \lambda_1 - m_3c_3 = [3, -2, 1, 3] - [0, -1, 2, 0] = [3, -1, -1, 3], \\
   & s_2(\lambda_2) = \lambda_2 - m_2c_2 = [2, 1, -2, 2] - [-1, 2, -1, -1] = [3, -1, -1, 3].
  \end{split}
\end{equation*}

\begin{table}[H]
\footnotesize
\centering
\renewcommand{\arraystretch}{0.98}

\begin{tabular}{||c|c||c|c||c|c||c|c||}
\hline
  \multicolumn{2}{||c||}{}      &
  \multicolumn{2}{c||}{}      &   \multicolumn{2}{c||}{}       & \multicolumn{2}{c||}{}      \\
  \multicolumn{2}{||c||}{$B_7$} &
  \multicolumn{2}{c||}{$D_8$}   &   \multicolumn{2}{c||}{$E_7$}  & \multicolumn{2}{c||}{$B_8$}      \\
  \multicolumn{2}{||c||}{}      &
  \multicolumn{2}{ c||}{}      &   \multicolumn{2}{c||}{}       & \multicolumn{2}{c||}{}      \\
  \hline
 level  & size  & level  & size  & level & size  & level & size  \\
\hline
  0, 49 &     1 & 0,  56 &      1  &  0, 63 &     1  & 0,  64 &      1 \\
  1, 48 &     7 & 1,  55 &      8  &  1, 62 &     7  & 1,  63 &      8 \\
  2, 47 &    27 & 2,  54 &     35  &  2, 61 &    27  & 2,  62 &     35 \\
  3, 46 &    77 & 3,  53 &    112  &  3, 60 &    77  & 3,  61 &    112 \\
  4, 45 &   181 & 4,  52 &    293  &  4, 59 &   182  & 4,  60 &    293 \\
  5, 44 &   371 & 5,  51 &    664  &  5, 58 &   378  & 5,  59 &    664 \\
  6, 43 &   686 & 6,  50 &   1350  &  6, 57 &   713  & 6,  58 &   1350 \\
  7, 42 &  1170 & 7,  49 &   2520  &  7, 56 &  1247  & 7,  57 &   2520 \\
  8, 41 &  1869 & 8,  48 &   4388  &  8, 55 &  2051  & 8,  56 &   4389 \\
  9, 40 &  2827 & 9,  47 &   7208  &  9, 54 &  3205  & 9,  55 &   7216 \\
 10, 39 &  4082 & 10, 46 &  11263  & 10, 53 &  4975  & 10, 54 &  11298 \\
 11, 38 &  5662 & 11, 45 &  16848  & 11, 52 &  6909  & 11, 53 &  16960 \\
 12, 37 &  7581 & 12, 44 &  24248  & 12, 51 &  9632  & 12, 52 &  24541 \\
 13, 36 &  9835 & 13, 43 &  33712  & 13, 50 & 13040  & 13, 51 &  34376 \\
 14, 35 & 12399 & 14, 42 &  45425  & 14, 49 & 17194  & 14, 50 &  46775 \\
 15, 34 & 15225 & 15, 41 &  59480  & 15, 48 & 22134  & 15, 49 &  62000 \\
 16, 33 & 18242 & 16, 40 &  75853  & 16, 47 & 27874  & 16, 48 &  80241 \\
 17, 32 & 21358 & 17, 39 &  94384  & 17, 46 & 34398  & 17, 47 & 101592 \\
 18, 31 & 24464 & 18, 38 & 114766  & 18, 45 & 41657  & 18, 46 & 126029 \\
 19, 30 & 27440 & 19, 37 & 136544  & 19, 44 & 49567  & 19, 45 & 153392 \\
 20, 29 & 30162 & 20, 36 & 159125  & 20, 43 & 58009  & 20, 44 & 183373 \\
 21, 28 & 32150 & 21, 35 & 181800  & 21, 42 & 66831  & 21, 43 & 215512 \\
 22, 27 & 34376 & 22, 34 & 203777  & 22, 41 & 75852  & 22, 42 & 249201 \\
 23, 26 & 35672 & 23, 33 & 224224  & 23, 40 & 84868  & 23, 41 & 283704 \\
 24, 25 & 36336 & 24, 32 & 242318  & 24, 39 & 93659  & 24, 40 & 318171 \\
        &       & 25, 31 & 257295  & 25, 38 & 101997 & 25, 39 & 351680 \\
        &       & 26, 30 & 268504  & 26, 37 & 109655 & 26, 38 & 383270 \\
        &       & 27, 29 & 275440  & 27, 36 & 116417 & 27, 37 & 411984 \\
        &       &  28   & 277788  & 28, 35 & 122087 & 28, 36 & 436913 \\
        &       &       &         & 29, 34 & 126497 & 29, 35 & 457240 \\
        &       &       &         & 30, 33 & 129514 & 30, 34 & 472281 \\
        &       &       &         & 31, 32 & 131046 & 31, 33 & 481520 \\
        &       &       &         &        &        &   32   & 484636 \\
\hline
\multicolumn{2}{||c||}{total 645120} &  \multicolumn{2}{c||}{total 5169960} &
\multicolumn{2}{c||}{total 2903040 } & \multicolumn{2}{c||}{total 10321920}      \\
\multicolumn{2}{||c||}{time 59 sec} &  \multicolumn{2}{c||}{time 570 sec} &
\multicolumn{2}{c||}{time 269 sec}  & \multicolumn{2}{c||}{time 1153 sec}      \\
\hline
\end{tabular}
\vspace{2mm}\caption{The Weyl groups $B_7$, $D_8$, $E_7$, $B_8$:
level sizes and total runtime of the extended Snow's algorithm}
\label{tab_levels_BDE}
\end{table}

Therefore, $s_3(\lambda_1) =  s_2(\lambda_2)$ and $s_3s_2s_3(\lambda_0) = s_2s_3s_2(\lambda_0)
\footnotemark[1]$.
Thus, weight $[3, -1, -1, 3]$ can be obtained in different ways.
This means that both $s_3w_1 = s_3s_2s_3$ and $s_2w_2 = s_2s_3s_2$ must be added
to the list of level $3$, even if they are two different reduced expressions for the same element.
\footnotetext[1]{The last relation also follows from the well-known braid relation
$s_3s_2s_3 = s_2s_3s_2$.}
\begin{figure}[h]
\centering
\includegraphics[scale=0.27]{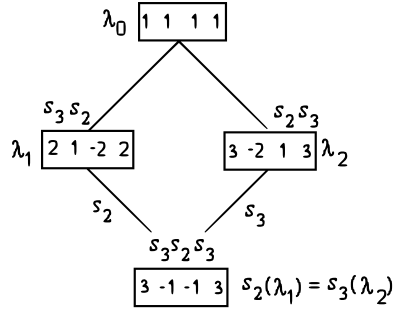}
\caption{Two different reduced expressions which are equal.}
\label{Bo_D4}
\end{figure}

This is an example of the repetition problem.
Snow's algorithm solves this problem with the following statement.
\begin{theorem}[Snow, \cite{Sn90}]
  \label{th_snow}
  Let $L_k$ be the $k$th level in the orbit $W\cdot\mu$ of a dominant weight $\mu \in \overline{C}$.
  Then, for each $\xi = (x_1,\dots, x_l) \in L_{k+1}$, there exists a unique $\nu \in L_k$ and
  a unique simple reflection $s_i$ such that $s_i(\nu) = \xi$ and $x_i \geq 0$ for  $j > i$.
  In particular, the next level $L_{k+1}$  can be constructed without repetitions from the weights
  $\nu \in L_k$ by adding $s_i(\nu)$ to $L_{k+1}$ if and only if the $i$th coordinate of $\nu$ is positive
  and the coordinates of $s_i(\nu)$ after the $i$th are nonnegative:
\begin{equation}
  \begin{split}
    L_{k+1} & = \{ s_i(\nu) = (x_1, \dots, x_l \mid i = 1,\dots, l, \\
            &  \nu = (y_1, \dots, y_l) \in L_k, y_i > 0, x_k \geq 0, j > i  \}.
  \end{split}
\end{equation}
\end{theorem}
~\\

\subsubsection{Application of Theorem \ref{th_snow} to example \ref{ex_repet_problem}}
Here $\xi = [3, -1, -1, 3]$.
For $v = \lambda_1$ and reflection $s_3$, we have $i = 3$ and $x_4 > 0$.
By Theorem \ref{th_snow} the element $s_3s_2s_3$ is added to level 3, see Table \ref{tab_level_3_1},
element 10.
On the other hand, for $v = \lambda_2$ and reflection $s_2$, we have $i = 2$ and
$x_3 < 0$. Then, the element $s_2s_3s_2$, which is essentially another reduced expression
for $s_3s_2s_3$, is not added to level 3.

\section{\bf Extended Snow's algorithm: computation of inverse elements}
  \label{sec_snow_ext}

\subsection{Double identification}

Because the reduced expression is not unique,
we must use another element identification $w$ to recognize the inverse element.
The matrix of $w$ in the faithful representation can be chosen as such a requested identifier.
We store the following information ({\it class Element}) about each element $w$:

\begin{python}
        weight       ''' w(v), where v is dominant integral weight '''
        name         ''' reduced expression  w like s3.s4.s2.s3.s1 '''
        matr         ''' matrix  w, two-dimensional list '''
        name_inv     ''' reduced expression of inverse '''
        matr_inv     ''' inverse matrix '''
        n_in_lvl     ''' location of w in Level '''
        n_inv_in_lvl ''' location of inverse in Level '''
\end{python}
For a complete description of this class, see $\S$\ref{data_element}.
The pair ({\it name, matr}) forms the double identification of the element.
The question is why not use a weight which is just $1D$-array
rather than a matrix which is $2D$-array.  The reason is that
at the time of calculating the new element given the element $w$,
we do not know the weight of the inverse element $w^{-1}$.
However, we know the inverse matrix $w^{-1}$, and, at the same time,
we do not perform a very expensive matrix inversion procedure.
Let $i$ the index of the desired reflection in the list of reflections {\it refl}.
Then $refl[i]$ (resp. $s_i$) is the matrix (resp. the symbol) of this reflection.
All we have to do is
\begin{itemize}
\item multiply the given matrix $w$ on the left by $refl[i]$
and the inverse matrix $w^{-1}$ on the right  by the  same reflection.
\item add the symbol $s_i$ on the left to the reduced expression $w$,
and for the reduced expression $w^{-1}$ add the symbol $s_i$ on the right.
\end{itemize}
When implemented in Python, it looks like this:
\begin{python}
   new_name = 's' + str(i) + '.' + name
   new_matr =  np.mathmul(refl[i], matr)
   new_name_inv = name_inv + '.s' + str(i)
   new_matr = np.mathmul(matr_inv, refl[i])
\end{python}
see function {\it newElem} in \S\ref{py_levels}.
Here, {\it np.mathmul} is a function from the {\it Numpy} package for multiplying two matrices.
Dot '.' used as delimiter between
generators in string fields $name$, $name\_inv$ and $new\_name\_inv$.

\subsection{Dictionary whose key is a matrix}
The dictionary {\it dictElemsOfLevel} is used to exchange information between
any element $w$ and its inverse $w^{-1}$.  The dictionary key is the matrix from
{\it class Element}.
The matrix is presented as a two-dimensional list.
Since a list cannot be a dictionary key in Python, we are converting
the matrix to a string as follows:
\begin{python}
   key   = ''.join(str(i) for row in self.matr for i in row)
\end{python}
The dictionary value corresponding to this {\it key} is the location $n\_in\_lvl$
of the matrix in $level(\xi)$. See function {\it keyValAndKeyInv() } in $\S$\ref{data_element}.
Let $key$ (resp. $key\_inv$) be the key corresponding to the $new\_matr$ (resp. $new\_matr\_inv$).
In the calculation cycle new level $L_{k+1}$ by the level $L_k $, there are $3$ cases, see function
{\it findAllLevels\_to\_LvlK()} in $\S$\ref{py_levels}.
Each record of the dictionary is the pair ({\it key}, {\it value}), where {\it key} is the matrix
converted to string, and {\it value} is the location of $w$ in $L_{k+1}$.

\subsection{Exchange information between $w$ and $w^{-1}$}
The element $w$ leaves in the dictionary record about its location in $L_{k+1}$. The inverse
element $w^{-1}$ will read this record later. There are three typical cases:
~\\

{\it Case 1}. If the computed matrix {\it new\_matr} is of order $2$,
i.e., the matrix is inverse to itself, then no message should be left in the dictionary.
This is the simplest case. Here,
\begin{python}
   new_elm.n_inv_in_lvl = new_elm.n_in_lvl        #   (1)
\end{python}

{\it Case 2}:
Suppose, after checking the key of the element $w$, it turned out that
{\it the key is not in the dictionary}. This means that the inverse element
will appear later in the calculation loop.
Then, the record about the location of $w$ is recorded in the dictionary.
\begin{python}
   dictElemsOfLevel[key_inv] = n_in_lvl           #   (2)
\end{python}
 The inverse element $w^{-1}$ will read this record later, see (3).
~\\

{\it Case 3}:
Suppose \underline{the key is in the dictionary}.
This means that the inverse element left an exact record about its location, see (2):
\begin{python}
   n_in_lvl = dictElemsOfLevel[key]               #   (3)
\end{python}
Then, there is no need to write any information to the  dictionary
because both $new\_elem$ and $new\_elem\_inv$ are already informed about
each other's location:
\begin{python}
   new_elem_inv = new_level[n_in_lvl]
   new_elem_inv.n_inv_in_lvl = new_elm.n_in_lvl
   new_elm.n_inv_in_lvl = new_elem_inv.n_in_lvl
\end{python}

\begin{remark}{\rm
Note that the keys will be recorded into the dictionary only for {\it Case 2}.
Let $\nu$ be number of records of some level $L_k$,
let $\omega_2$ be the number of elements of order $2$ in $L_k$,
and $\delta$ be the number of elements of $L_k$ in the dictionary at the end of
the run cycle. Then
\begin{equation}
     \delta = (\nu - \omega_2)/2.
\end{equation}
The number of elements of any level in the dictionary
will always be less than half of all elements of this level.
}
\end{remark}

\subsection{Complexity}

The extended Snow's algorithm makes it possible to find pairs of
inverse elements in the same work cycle. Note that time-complexity of Snow's algorithm  is O(n),
and the extended Snow's algorithm also has the time-complexity O(n), there $n$ is the order
of the Weyl group.

The algorithm uses the dictionary mechanism which is very effective in Python.
The very frequently dictionary operation keys() has the time-complexity O(n) in Python-2  and
O(1) in Python-3. The average time-complexity for operations get() and set() are also O(1), \cite{Py22}.

~\\
\section{\bf Conjugacy classes in $W(D_4)$}
In this section, we consider different representations of the conjugacy classes in $W(D_4)$.

\subsection{Conjugacy classes of $W(D_4)$ represented by Carter diagrams}

First, we will see why in Table \ref{tab_of_all_ccls} the representative element
\begin{equation}
 \label{repr_ccl_6}
  s_1s_2s_3s_4s_2s_1s_2s_3s_4s_2s_3s_4
\end{equation}
of the conjugacy class $12$ is represented as $4$ unconnected vertices (root subset $4A_1$),
and the representative element
\begin{equation}
 \label{repr_ccl_11}
 s_3s_2s_4s_3s_2s_1
\end{equation}
of the conjugacy class $11$ is represented by the Carter diagram $D_4(a_1)$.

For more convenient work with roots of the root system $D_4$, we change the notation
of vertices from $i$ to $\alpha_i$.
We use the Bourbaki numbering of the vertices of the Dynkin diagram $D_4$:
The reflection $s_{\alpha_2}$ does not commute with reflections $s_{\alpha_i}$, $i = 1,3,4$,
while the reflections $s_{\alpha_1}$, $s_{\alpha_3}$, $s_{\alpha_4}$ commute with each other,
see Carter diagram in Table \ref{tab_of_all_ccls},line 11.

\begin{table}[H]
\footnotesize
\centering
\renewcommand{\arraystretch}{0.9}
   \begin{tabular}{|c|c|c|c|c|c|c|}
     \hline
        $N^{\circ}$  &  Carter                  & Representative
                                 & Elms  & Root   & Order  &  Signed  \\
                     &  diagram\footnotemark[1] &   element
                                 &       & subset &        & cycle-type\footnotemark[2]\\
     \hline
        0     &
         -       &  $e$         &  1        & $\emptyset$ & 1  & [1111] \\
     \hline
        1     &
       $\begin{array}{c}
       \includegraphics[scale=0.4]{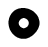}
       \end{array}$
                &  $s_1$         &  12        & $A_1$  & 2  & [211]  \\
     \hline
        2     &
       $\begin{array}{c}
       \includegraphics[scale=0.4]{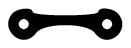}
       \end{array}$
                 &  $s_1s_2$      &  32        & $A_2$ & 3 & [31] \\
     \hline
        3     &
       $\begin{array}{c}
       \includegraphics[scale=0.4]{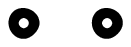}
       \end{array}$
                &  $s_1s_3$      &  6         & $2A_1$ & 2 & [22] \\
     \hline
        4     &
       $\begin{array}{c}
       \includegraphics[scale=0.4]{D2.eps}
       \end{array}$
         &  $s_1s_4$      &  6         & $2A_1$ & 2 &  [22] \\
     \hline
        5     &
       $\begin{array}{c}
       \includegraphics[scale=0.4]{D2.eps}
       \end{array}$
        &  $s_3s_4$      &  6         & $D_2$ & 2 & [$\bar{1}\bar{1}$11] \\
     \hline
        6     &
       $\begin{array}{c}
       \includegraphics[scale=0.4]{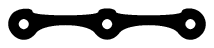}
       \end{array}$
              &  $s_1s_2s_3$   &  24    & $A_3$ &  4 & [4] \\
     \hline
        7     &
       $\begin{array}{c}
       \includegraphics[scale=0.4]{A3.eps}
       \end{array}$
               &  $s_1s_2s_4$   &  24    & $A_3$ & 4 & [4] \\
     \hline
        8     &
       $\begin{array}{c}
       \includegraphics[scale=0.4]{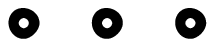}
       \end{array}$
        &  $s_1s_3s_4$   &  12        & $3A_1$  & 2 &  [2$\bar{1}\bar{1}$] \\
     \hline
        9    &
       $\begin{array}{c}
       \includegraphics[scale=0.4]{A3.eps}
       \end{array}$
         &  $s_3s_2s_4$   &  24        & $D_3$ & 4 & [$\bar{2}\bar{1}$1] \\
     \hline
        10    &
       $\begin{array}{c}
       \includegraphics[scale=0.3]{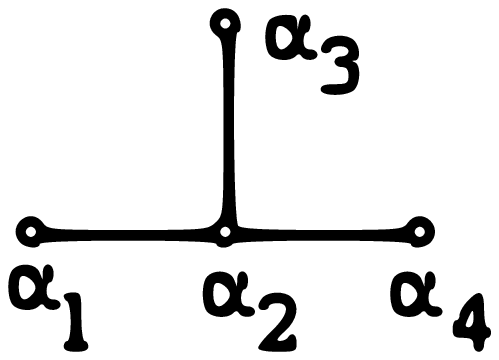}
       \end{array}$
         &  $s_1s_4s_2s_3$   &  32     & $D_4$ & 6 & [$\bar{3}\bar{1}$] \\
     \hline
        11    &
       $\begin{array}{c}
       \includegraphics[scale=0.3]{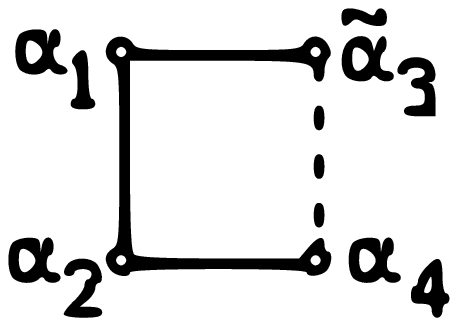}
       \end{array}$
                &  $s_3s_2s_4s_3s_2s_1$   &  12  & $D_4(a_1)$ & 4 & [$\bar{2}\bar{2}$]\\
     \hline
        12     &
       $\begin{array}{c}
       \includegraphics[scale=0.4]{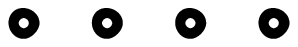}
       \end{array}$
        &  $s_1s_2s_3s_4s_2s_1s_2s_3s_4s_2s_3s_4$  & 1 & $4A_1$ & 2 &
             [$\bar{1}\bar{1}\bar{1}\bar{1}$]\\
     \hline
 \end{tabular}
 \vspace{2mm}
\caption{Conjugacy classes in the Weyl groups $D_4$, see Tables \ref{tableCCL_0123} - \ref{tableCCL_12}}
\label{tab_of_all_ccls}
\centering
   \begin{tabular}{|c|c|c|c|c|c|}
     \hline
   Order          &  $1$         & $2$  & $3$  & $4$  & $6$  \\
         \hline
   Elements      &  $1$         & $43$ & $32$ & $84$ & $32$ \\
         \hline
 \end{tabular}
  \vspace{1mm}
\caption{Weyl groups $D_4$. Partitioning by element orders}
\label{tab_part_by_order}
\end{table}
\footnotetext[1]{For an explanation of the Carter diagram $D_4(a_1)$
with a dotted edge (in the 11th conjugacy class), see \cite[$\S$1.1.1]{St17}.}
\footnotetext[2]{For the definition of signed cycle-types, see \cite[$\S$7]{Ca72}.}
~\\
For any pair of non-orthogonal roots $\alpha$ and $\beta$,
such that $(\alpha, \beta) = -1$, the following relations hold:
\begin{equation}
 \label{pair_non_orthog_roots}
 \begin{split}
   & s_{\beta}s_{\alpha}s_{\beta} = s_{s_\beta({\alpha})} = s_{\alpha+\beta},
   \text{ and } s_{\beta}s_{\alpha} = s_{\alpha+\beta}s_{\beta}, ~~
    s_{\alpha}s_{\beta} = s_{\beta}s_{\alpha+\beta}.  \\
   & (s_{\beta}s_{\alpha})^3 = 1,  \text{ since }
   (s_{\beta}s_{\alpha})^3 =
   (s_{\beta}s_{\alpha}s_{\beta})(s_{\alpha}s_{\beta}s_{\alpha}) =
   s_{\alpha+\beta}^2 = 1.
 \end{split}
\end{equation}

\subsubsection{Conjugacy class $11$, Carter diagram $D_4(a_1)$} The representative element
$w =  s_3s_2s_4s_3s_2s_1$
is the first element of conjugacy class $11$, see Table \ref{tableCCL_11}.
Using the roots from the root system as indices, we get the following expression for $w$:
\begin{equation*}
   w =  s_{\alpha_3}s_{\alpha_2}s_{\alpha_4}s_{\alpha_3}s_{\alpha_2}s_{\alpha_1} =
        s_{\alpha_3}s_{\alpha_2}s_{\alpha_3}s_{\alpha_4}s_{\alpha_2}s_{\alpha_1} =
        s_{\alpha_2 + \alpha_3}s_{\alpha_4}s_{\alpha_2}s_{\alpha_1}.
\end{equation*}
Further,
\begin{equation}
  \label{elem_w}
  w \stackrel{s_{\alpha_1}}{\simeq} s_{\alpha_1}s_{\alpha_2 + \alpha_3}s_{\alpha_4}s_{\alpha_2},
\end{equation}
where, the notaion $\stackrel{A}{\simeq}$ means conjugacy by the element $A$.
The element \eqref{elem_w} can be transformed as follows:
\begin{equation}
  \label{D4a1_4_roots}
   w = s_{\alpha_1}s_{\alpha_2 + \alpha_3}s_{\alpha_4}s_{\alpha_2} =
    s_{\alpha_2 + \alpha_3 + \alpha_1}s_{\alpha_1}s_{\alpha_4}s_{\alpha_2} = \\
         s_{\widetilde\alpha_3}s_{\alpha_1}s_{\alpha_4}s_{\alpha_2},
\end{equation}
where $\widetilde\alpha_3 = -(\alpha_1 + \alpha_2 + \alpha_3)$.
The element $w$ is represented by the Carter diagram $D_4(a_1)$, where
the dotted edge $\{\widetilde\alpha_3, \alpha_4\}$ corresponds to the inner product
$(\widetilde\alpha_3, \alpha_4) = 1$, see \cite{St17}.

\subsubsection{Conjugacy class $12$, four unconnected vertices}
The element \eqref{repr_ccl_6} looks like this:
\begin{equation*}
w =  s_{\alpha_1}s_{\alpha_2}s_{\alpha_3}s_{\alpha_4}s_{\alpha_2}s_{\alpha_1}s_{\alpha_2}s_{\alpha_3}s_{\alpha_4}s_{\alpha_2}s_{\alpha_3}s_{\alpha_4}.
\end{equation*}
First of all, according to \eqref{pair_non_orthog_roots}, we change
$s_{\alpha_2}s_{\alpha_1}s_{\alpha_2}$ to $s_{\alpha_2 + \alpha_1}$
and $s_{\alpha_4}s_{\alpha_2}s_{\alpha_4}$ to $s_{\alpha_2 + \alpha_4}$.  Then
\begin{equation*}
w =  s_{\alpha_1}s_{\alpha_2}s_{\alpha_3}s_{\alpha_4}s_{\alpha_2 + \alpha_1}s_{\alpha_3}s_{\alpha_2 + \alpha_4}s_{\alpha_3}.
\end{equation*}
Further, by \eqref{pair_non_orthog_roots},
we change $s_{\alpha_3}s_{\alpha_2 + \alpha_4}s_{\alpha_3}$ to
$s_{\alpha_2 + \alpha_4 + \alpha_3}$, and $s_{\alpha_4}s_{\alpha_2 + \alpha_1}$ to
$s_{\alpha_2 + \alpha_1}s_{\alpha_2 + \alpha_1 + \alpha_4}$. Thus,
\begin{equation*}
w =  s_{\alpha_1}s_{\alpha_2}s_{\alpha_3}s_{\alpha_2 + \alpha_1}s_{\alpha_2 + \alpha_1 + \alpha_4}s_{\alpha_2 + \alpha_4 + \alpha_3}.
\end{equation*}
Similarly, we replace $s_{\alpha_3}s_{\alpha_2 + \alpha_1}$  with
$s_{\alpha_2 + \alpha_1}s_{\alpha_2 + \alpha_1 + \alpha_3}$, we get
\begin{equation*}
    w =  s_{\alpha_1}s_{\alpha_2}s_{\alpha_2 + \alpha_1}
       s_{\alpha_2 + \alpha_1 + \alpha_3}s_{\alpha_2 + \alpha_1 + \alpha_4}s_{\alpha_2 + \alpha_4 + \alpha_3}.
\end{equation*}
At last, since $s_{\alpha_2}s_{\alpha_2 + \alpha_1} = s_{\alpha_1}s_{\alpha_2}$, we have
$s_{\alpha_1}s_{\alpha_2}s_{\alpha_2 + \alpha_1} = s_{\alpha_2}$ and
\begin{equation}
 \label{4_roots}
w =  s_{\alpha_2}s_{\alpha_2 + \alpha_1 + \alpha_3}s_{\alpha_2 + \alpha_1 + \alpha_4}s_{\alpha_2 + \alpha_4 + \alpha_3}.
\end{equation}
Note that in eq. \eqref{4_roots}, there are four  mutually orthogonal roots:
\begin{equation}
 \label{4_roots_2}
\alpha_2, \quad \alpha_2 + \alpha_1 + \alpha_3, \quad
\alpha_2 + \alpha_1 + \alpha_4, \quad \alpha_2 + \alpha_4 + \alpha_3
\end{equation}
The subset \eqref{4_roots_2} is represented by $4$ unconnected vertices, i.e., $4A_1$.

\subsection{Conjugacy classes of $W(D_4)$ represented by signed cycle-types}
In this section, we consider the representation of conjugacy classes $8$--$12$ of  Table \ref{tab_of_all_ccls}
using the signed cycle-types. According to Bourbaki's notaion:
\begin{equation*}
   s_{\alpha_1} = s_{e_1-e_2}, \; s_{\alpha_2} = s_{e_2-e_3}, \; s_{\alpha_3} = s_{e_3-e_4}, \;  s_{\alpha_4} = s_{e_3+e_4}.
\end{equation*}
We will use the following mappings:
\begin{equation}
  \label{formulas_ei}
\footnotesize
 \begin{array}{ccc}
   s_{e_i-e_j}:
    \begin{cases}
     e_i \longmapsto e_j \\
     e_j \longmapsto e_i \\
    \end{cases} \;
   s_{e_i+e_j}:
    \begin{cases}
     e_i \longmapsto -e_j \\
     e_j \longmapsto -e_i \\
    \end{cases} \;
   s_{e_i-e_j}s_{e_i+e_j}:
    \begin{cases}
     e_i \longmapsto -e_i \\
     e_j \longmapsto -e_j \\
    \end{cases},
 \end{array}
\end{equation}
see \cite[Ch.VI, $\S4$, $n^\circ$8]{Bo02}.

\subsubsection{Conjugacy class $8$, signed cycle-type $[2\bar{1}\bar{1}]$}
Consider representative element
$s_{\alpha_1}s_{\alpha_3}s_{\alpha_4}$. Let us find the signed cycle-type of this element.
By \eqref{formulas_ei},
$s_{e_1-e_2}$ permutes $e_1$ and $e_2$, i.e., $s_{e_1-e_2}$  acts as permutation $(12)$.
Further, the product $s_{e_3-e_4}s_{e_3+e_4}$ maps $e_3$ to $-e_3$ and $e_4$ to $-e_4$, t.e., acts as
the pair of negative cycles [$\bar{1}\bar{1}$]. All together gives [2$\bar{1}\bar{1}$].

\subsubsection{Conjugacy class $9$, signed cycle-type $[\bar{2}\bar{1}1]$}
Here, the representative element is $s_{\alpha_3}s_{\alpha_2}s_{\alpha_4}$. By \eqref{formulas_ei}
$s_{\alpha_2}$ permutes $e_2$ and $e_3$; $s_{\alpha_3}$ permutes $e_3$ and $e_4$. At last,
$s_{\alpha_4}$ maps $e_4$ to $-e_3$ and $e_3$ to $-e_4$. Then,
\begin{equation*}
     s_{\alpha_3}s_{\alpha_2}s_{\alpha_4}:
     \begin{cases}
        e_2 \longmapsto e_4, \\
        e_3 \longmapsto -e_3, \\
        e_4 \longmapsto -e_2.
     \end{cases}
\end{equation*}
The second mapping corresponds to the negative cycle [$\bar{1}$].
The first and third mappings form the cycle $e_2 \longmapsto e_4 \longmapsto -e_2$, i.e.,
the negative cycle [$\bar{2}$].  Thus, we get the signed cycle-type [$\bar{2}\bar{1}$],
or, that is the same, [$\bar{2}\bar{1}$1].
By \cite[Prop. 25]{Ca72}, [$\bar{i}\bar{1}$] corresponds to the Carter diagram $D_{i+1}$.
In our case, we get $D_3$.

\subsubsection{Conjugacy class $10$, signed cycle-type $[\bar{3}\bar{1}]$}
The representative element
\begin{equation}
 s_{\alpha_1}s_{\alpha_4}s_{\alpha_2}s_{\alpha_3} = s_{e_1-e_2}s_{e_3+e_4}s_{e_2-e_3}s_{e_3-e_4}
\end{equation}
acts as follows:
\begin{equation*}
     s_{\alpha_1}s_{\alpha_4}s_{\alpha_2}s_{\alpha_3}:
     \begin{cases}
        e_1 \longmapsto e_2, \quad e_3 \longmapsto -e_3, \\
        e_2 \longmapsto -e_4, \quad e_4 \longmapsto e_1. \\
     \end{cases}
\end{equation*}
The mapping $e_3 \longmapsto -e_3$ corresponds to the negative cycle [$\bar{1}$].
The remaining mappings
form the cycle $e_1 \longmapsto e_2 \longmapsto -e_4 \longmapsto -e_1$, i.e., the
negative cycle [$\bar{3}$]. So, we get the signed cycle-type [$\bar{3}\bar{1}$].
As above, by \cite[Prop. 25]{Ca72}, the signed cycle-type [$\bar{3}\bar{1}$] corresponds to $D_4$.

\subsubsection{Conjugacy class $11$, signed cycle-type $[\bar{2}\bar{2}]$}
By \eqref{D4a1_4_roots}, the representative element
\begin{equation*}
   s_{\alpha_2 + \alpha_3 + \alpha_1}s_{\alpha_1}s_{\alpha_4}s_{\alpha_2} =
   s_{e_1 - e_4}s_{e_1 - e_2}s_{e_3 + e_4}e_{e_2 - e_3}
\end{equation*}
acts as follows:
\begin{equation*}
   s_{\alpha_2 + \alpha_3 + \alpha_1}s_{\alpha_1}s_{\alpha_4}s_{\alpha_2}:
    \begin{cases}
      e_1 \longmapsto e_2, \\
      e_2 \longmapsto -e_4 \longmapsto -e_1, \\
      e_3 \longmapsto e_4, \\
      e_4 \longmapsto -e_3.
    \end{cases}
\end{equation*}
The first and second mappings
form the cycle $e_1 \longmapsto e_2 \longmapsto -e_1$, i.e., the
negative cycle [$\bar{2}$]. The third and fourth mappings
form the cycle $e_3 \longmapsto e_4 \longmapsto -e_3$, which is also the
negative cycle [$\bar{2}$]. Thus, we get the signed cycle-type [$\bar{2}\bar{2}$].

\subsubsection{Conjugacy class $12$, signed cycle-type $[\bar{1}\bar{1}\bar{1}\bar{1}]$}
By \eqref{4_roots} the representative element is as follows
\begin{equation*}
 s_{\alpha_2}s_{\alpha_2 + \alpha_1 + \alpha_3}
    s_{\alpha_2 + \alpha_1 + \alpha_4}s_{\alpha_2 + \alpha_4 + \alpha_3} =
 e_{e_2 - e_3}e_{e_1 - e_4}e_{e_1 + e_4}e_{e_2 + e_3}.
\end{equation*}
Since $s_{e_i - e_j}s_{e_i + e_j}$ maps $e_i$ to $-e_i$ and $e_j$ to $-e_j$,
we get
\begin{equation*}
    s_{\alpha_2}s_{\alpha_2 + \alpha_1 + \alpha_3}
    s_{\alpha_2 + \alpha_1 + \alpha_4}s_{\alpha_2 + \alpha_4 + \alpha_3}: e_i \longmapsto -e_i
    \text{ for } i = 1,2,3,4.
\end{equation*}
This corresponds to the signed cycle-type $[\bar{1}\bar{1}\bar{1}\bar{1}]$.


\section{\bf Weyl group $D_4$. Partitioning by levels }
   \label{sec_levels_D4}

\begin{table}[H]
  \centering
  \footnotesize
  \renewcommand{\arraystretch}{0.7}
\begin{tabular} {||c|c|c|c||}
    \hline
     \footnotesize $N^o$ & \footnotesize Weight & \footnotesize Element & \footnotesize Matrix \\
    \hline
     \footnotesize $0$ & \footnotesize $1,1,1,1$ & $e$
                                    & \footnotesize
                                   $\begin{bmatrix*}[r]
                                     1 & 0 & 0 & 0 \\
                                     0 & 1 & 0 & 0 \\
                                     0 & 0 & 1 & 0 \\
                                     0 & 0 & 0 & 1 \\
                                   \end{bmatrix*}$ \\
    \hline
  \end{tabular}
  \caption{\footnotesize Weyl group $D_4$, level 0, element $e$}

  \begin{tabular} {||c|c|c|c||c|c|c|c||}
    \hline
     $N^o$ & Weight & Element & Matrix & $N^o$ & Weight & Element & Matrix \\
    \hline
   \footnotesize $0$ & \footnotesize  $-1,2,1,1$ & $s_1$
      & \footnotesize $\begin{bmatrix*}[r]
          -1 & 1 & 0 & 0 \\
          0  & 1 & 0 & 0 \\
          0  & 0 & 1 & 0 \\
          0  & 0 & 0 & 1 \\
         \end{bmatrix*}$ &
   \footnotesize $1$ & \footnotesize $2,-1,2,2$ & $s_2$
      &  \footnotesize $\begin{bmatrix*}[r]
          1 &  0 & 0 & 0 \\
          1 & -1 & 1 & 1 \\
          0 &  0 & 1 & 0 \\
          0 &  0 & 0 & 1 \\
         \end{bmatrix*}$ \\
    \hline
   \footnotesize $2$ & \footnotesize $1,2,-1,1$ & $s_3$
      & \footnotesize $\begin{bmatrix*}[r]
          1 & 0 &  0 & 0 \\
          0 & 1 &  0 & 0 \\
          0 & 1 & -1 & 0 \\
          0 & 0 &  0 & 1  \\
         \end{bmatrix*}$ &
    \footnotesize $3$ & \footnotesize $1,2,1,-1$ & $s_4$
      &  \footnotesize $\begin{bmatrix*}[r]
          1 &  0 & 0 & 0 \\
          0 &  1 & 0 & 0 \\
          0 &  0 & 1 & 0 \\
          0 &  1 & 0 & -1 \\
         \end{bmatrix*}$ \\
    &&&&&&&\\
    \hline
  \end{tabular}
  \caption{\footnotesize Weyl group $D_4$, level 1, elements 0-3}

  \begin{tabular} {||c|c|c|c||c|c|c|c||}
    \hline
    &&&&&&&\\
     $N^o$ & Weight & Element  & Matrix & $N^o$ & Weight & Element   & Matrix \\
           &        & \& order &        &       &        & \& order  &     \\
    &&&&&&&\\
    \hline
    &&&&&&&\\
    \footnotesize $\begin{array}{c}
        3 \\
        \\
        (0) \footnotemark[1]
     \end{array}$
      & \footnotesize $-2,1,2,2$
      & $\begin{array}{c}
           s_1s_2 \\
           \\
           3 \\
         \end{array}$
      & \footnotesize $\begin{bmatrix*}[r]
          0 & -1 & 1 & 1 \\
          1 & -1 & 1 & 1 \\
          0 &  0 & 1 & 0 \\
          0 &  0 & 0 & 1 \\
         \end{bmatrix*}$ &
      \footnotesize $\begin{array}{c}
        6 \\
        \\
        (4)
     \end{array}$
      & \footnotesize $3,-2,1,3$
      & $\begin{array}{c}
           s_2s_3 \\
           \\
           3 \\
         \end{array}$
      &  \footnotesize $\begin{bmatrix*}[r]
          1 & 0 & 0 & 0 \\
          1 & 0 & -1 & 1 \\
          0 & 1 & -1 & 0 \\
          0 & 0 & 0  & 1 \\
         \end{bmatrix*}$ \\
    &&&&&&&\\
    \hline
    &&&&&&&\\
\hline
  \end{tabular}
  \caption{\footnotesize Weyl group $D_4$, level 2, elements 0-8}
  \label{tab_levels_012}
\end{table}
\footnotetext[1]{Hereinafter, the number in this column  without parentheses (resp. in parentheses)
means the ordinal number of element (resp. inverse element) in
Tables \ref{tab_levels_012} - \ref{tab_level_11_12}. } 
\begin{table}[H]
\centering
\renewcommand{\arraystretch}{0.98}
     \begin{tabular}{||c|c|c|c|c||}
     \hline
          &&&& \\ 
    $N^o$  & Weight &  Element   & Matrix & Inverse   \\ 
           &        & \& order  &        &           \\ 
     \hline
     \footnotesize
     $\begin{array}{c}
         0 \\ 
          \\ 
         (9) \\ 
     \end{array}$
&     \footnotesize
 $[1, 1, -3, 3]$&     $\begin{array}{c}
        s_3s_2s_1 \\ 
             \end{array}$
&     \footnotesize
     $\begin{bmatrix*}[r]
        -1 & 1 & 0 & 0 \\ 
        -1 & 0 & 1 & 1 \\ 
        -1 & 0 & 0 & 1 \\ 
        0 & 0 & 0 & 1 \\ 
     \end{bmatrix*}$
&     $\begin{array}{c}
        s_1s_2s_3 \\ 
     \end{array}$
   \\ 
     \hline
     \footnotesize
     $\begin{array}{c}
         1 \\ 
          \\ 
         (13) \\ 
     \end{array}$
&     \footnotesize
 $[1, 1, 3, -3]$&     $\begin{array}{c}
        s_4s_2s_1 \\ 
             \end{array}$
&     \footnotesize
     $\begin{bmatrix*}[r]
        -1 & 1 & 0 & 0 \\ 
        -1 & 0 & 1 & 1 \\ 
        0 & 0 & 1 & 0 \\ 
        -1 & 0 & 1 & 0 \\ 
     \end{bmatrix*}$
&     $\begin{array}{c}
        s_1s_2s_4 \\ 
     \end{array}$
   \\ 
     \hline
     \footnotesize
     $\begin{array}{c}
         2 \\ 
          \\ 
         (6) \\ 
     \end{array}$
&     \footnotesize
 $[2, -3, 2, 4]$&     $\begin{array}{c}
        s_2s_3s_1 \\ 
             \end{array}$
&     \footnotesize
     $\begin{bmatrix*}[r]
        -1 & 1 & 0 & 0 \\ 
        -1 & 1 & -1 & 1 \\ 
        0 & 1 & -1 & 0 \\ 
        0 & 0 & 0 & 1 \\ 
     \end{bmatrix*}$
&     $\begin{array}{c}
        s_3s_1s_2 \\ 
     \end{array}$
   \\
\hline
     \end{tabular}
\caption{\footnotesize Weyl group $D_4$, level 3, elements 0-13}
\label{tab_level_3_1}
\end{table}
\begin{table}[H]
\centering
\renewcommand{\arraystretch}{0.98}
     \begin{tabular}{||c|c|c|c|c||}
     \hline
          &&&& \\ 
    $N^o$  & Weight &  Element   & Matrix & Inverse   \\ 
           &        & \& order  &        &           \\ 
     \hline
     \footnotesize
     $\begin{array}{c}
         14 \\ 
          \\ 
         (11) \\ 
     \end{array}$
&     \footnotesize
 $[3, 1, -3, 1]$&     $\begin{array}{c}
        s_3s_2s_4 \\ 
             \end{array}$
&     \footnotesize
     $\begin{bmatrix*}[r]
        1 & 0 & 0 & 0 \\ 
        1 & 0 & 1 & -1 \\ 
        1 & 0 & 0 & -1 \\ 
        0 & 1 & 0 & -1 \\ 
     \end{bmatrix*}$
&     $\begin{array}{c}
        s_4s_2s_3 \\ 
     \end{array}$
   \\ 
     \hline
     \footnotesize
     $\begin{array}{c}
         15 \\ 
          \\ 
         (15) \\ 
     \end{array}$
&     \footnotesize
 $[3, -1, 3, -1]$&     $\begin{array}{c}
        s_4s_2s_4 \\ 
             \end{array}$
&     \footnotesize
     $\begin{bmatrix*}[r]
        1 & 0 & 0 & 0 \\ 
        1 & 0 & 1 & -1 \\ 
        0 & 0 & 1 & 0 \\ 
        1 & -1 & 1 & 0 \\ 
     \end{bmatrix*}$
&     $\begin{array}{c}
        s_4s_2s_4 \\ 
     \end{array}$
   \\ 
     \hline
     \end{tabular}
\caption{\footnotesize Weyl group $D_4$, level 3, elements 14-15}
\label{tab_level_3_2}
\end{table}

\begin{table}[H]
\centering
\renewcommand{\arraystretch}{0.98}
     \begin{tabular}{||c|c|c|c|c||}
     \hline
          &&&& \\ 
    $N^o$  & Weight &  Element   & Matrix & Inverse   \\ 
           &        & \& order  &        &           \\ 
     \hline
     \footnotesize
     $\begin{array}{c}
         0 \\ 
          \\ 
         (16) \\ 
     \end{array}$
&     \footnotesize
 $[1, 4, -3, -3]$&     $\begin{array}{c}
        s_4s_3s_2s_1 \\ 
             \end{array}$
&     \footnotesize
     $\begin{bmatrix*}[r]
        -1 & 1 & 0 & 0 \\ 
        -1 & 0 & 1 & 1 \\ 
        -1 & 0 & 0 & 1 \\ 
        -1 & 0 & 1 & 0 \\ 
     \end{bmatrix*}$
&     $\begin{array}{c}
        s_1s_2s_4s_3 \\ 
     \end{array}$
   \\ 
     \hline
     \footnotesize
     $\begin{array}{c}
         1 \\ 
          \\ 
         (13) \\ 
     \end{array}$
&     \footnotesize
 $[2, -1, -2, 4]$&     $\begin{array}{c}
        s_3s_2s_3s_1 \\ 
             \end{array}$
&     \footnotesize
     $\begin{bmatrix*}[r]
        -1 & 1 & 0 & 0 \\ 
        -1 & 1 & -1 & 1 \\ 
        -1 & 0 & 0 & 1 \\ 
        0 & 0 & 0 & 1 \\ 
     \end{bmatrix*}$
&     $\begin{array}{c}
        s_3s_1s_2s_3 \\ 
     \end{array}$
   \\ 
     \hline
     \footnotesize
     $\begin{array}{c}
         2 \\ 
          \\ 
         (20) \\ 
     \end{array}$
&     \footnotesize
 $[2, 1, 2, -4]$&     $\begin{array}{c}
        s_4s_2s_3s_1 \\ 
             \end{array}$
&     \footnotesize
     $\begin{bmatrix*}[r]
        -1 & 1 & 0 & 0 \\ 
        -1 & 1 & -1 & 1 \\ 
        0 & 1 & -1 & 0 \\ 
        -1 & 1 & -1 & 0 \\ 
     \end{bmatrix*}$
&     $\begin{array}{c}
        s_3s_1s_2s_4 \\ 
     \end{array}$
   \\
\hline
     \end{tabular}
\caption{\footnotesize Weyl group $D_4$, level 4, elements 0-13}
\label{tab_level_4_1}
\end{table}
\begin{table}[H]
\centering
\renewcommand{\arraystretch}{0.98}
     \begin{tabular}{||c|c|c|c|c||}
     \hline
          &&&& \\ 
    $N^o$  & Weight &  Element   & Matrix & Inverse   \\ 
           &        & \& order  &        &           \\ 
     \hline
     \footnotesize
     $\begin{array}{c}
         14 \\ 
          \\ 
         (4) \\ 
     \end{array}$
&     \footnotesize
 $[-3, 4, 1, -3]$&     $\begin{array}{c}
        s_4s_1s_2s_3 \\ 
             \end{array}$
&     \footnotesize
     $\begin{bmatrix*}[r]
        0 & 0 & -1 & 1 \\ 
        1 & 0 & -1 & 1 \\ 
        0 & 1 & -1 & 0 \\ 
        1 & 0 & -1 & 0 \\ 
     \end{bmatrix*}$
&     $\begin{array}{c}
        s_3s_2s_4s_1 \\ 
     \end{array}$
   \\ 
     \hline
     \footnotesize
     $\begin{array}{c}
         15 \\ 
          \\ 
         (17) \\ 
     \end{array}$
&     \footnotesize
 $[3, 2, -1, -3]$&     $\begin{array}{c}
        s_4s_3s_2s_3 \\ 
             \end{array}$
&     \footnotesize
     $\begin{bmatrix*}[r]
        1 & 0 & 0 & 0 \\ 
        1 & 0 & -1 & 1 \\ 
        1 & -1 & 0 & 1 \\ 
        1 & 0 & -1 & 0 \\ 
     \end{bmatrix*}$
&     $\begin{array}{c}
        s_3s_2s_4s_3 \\ 
     \end{array}$
   \\ 
     \hline
     \footnotesize
     $\begin{array}{c}
         16 \\ 
          \\ 
         (0) \\ 
     \end{array}$
&     \footnotesize
 $[-4, 1, 2, 2]$&     $\begin{array}{c}
        s_1s_2s_4s_3 \\ 
             \end{array}$
&     \footnotesize
     $\begin{bmatrix*}[r]
        0 & 1 & -1 & -1 \\ 
        1 & 1 & -1 & -1 \\ 
        0 & 1 & -1 & 0 \\ 
        0 & 1 & 0 & -1 \\ 
     \end{bmatrix*}$
&     $\begin{array}{c}
        s_4s_3s_2s_1 \\ 
     \end{array}$
   \\
\hline
     \end{tabular}
\caption{\footnotesize Weyl group $D_4$, level 4, elements 14-22}
\label{tab_level_4_2}
\end{table}

\begin{table}[H]
\centering
\renewcommand{\arraystretch}{0.98}
     \begin{tabular}{||c|c|c|c|c||}
     \hline
          &&&& \\ 
    $N^o$  & Weight &  Element   & Matrix & Inverse   \\ 
           &        & \& order  &        &           \\ 
     \hline
     \footnotesize
     $\begin{array}{c}
         0 \\ 
          \\ 
         (11) \\ 
     \end{array}$
&     \footnotesize
 $[5, -4, 1, 1]$&     $\begin{array}{c}
        s_2s_4s_3s_2s_1 \\ 
             \end{array}$
&     \footnotesize
     $\begin{bmatrix*}[r]
        -1 & 1 & 0 & 0 \\ 
        -2 & 1 & 0 & 0 \\ 
        -1 & 0 & 0 & 1 \\ 
        -1 & 0 & 1 & 0 \\ 
     \end{bmatrix*}$
&     $\begin{array}{c}
        s_1s_2s_4s_3s_2 \\ 
     \end{array}$
   \\ 
     \hline
     \footnotesize
     $\begin{array}{c}
         1 \\ 
          \\ 
         (20) \\ 
     \end{array}$
&     \footnotesize
 $[2, 3, -2, -4]$&     $\begin{array}{c}
        s_4s_3s_2s_3s_1 \\ 
             \end{array}$
&     \footnotesize
     $\begin{bmatrix*}[r]
        -1 & 1 & 0 & 0 \\ 
        -1 & 1 & -1 & 1 \\ 
        -1 & 0 & 0 & 1 \\ 
        -1 & 1 & -1 & 0 \\ 
     \end{bmatrix*}$
&     $\begin{array}{c}
        s_3s_1s_2s_4s_3 \\ 
     \end{array}$
   \\ 
     \hline
     \footnotesize
     $\begin{array}{c}
         2 \\ 
          \\ 
         (17) \\ 
     \end{array}$
&     \footnotesize
 $[3, -1, -3, 3]$&     $\begin{array}{c}
        s_3s_2s_4s_3s_1 \\ 
             \end{array}$
&     \footnotesize
     $\begin{bmatrix*}[r]
        -1 & 1 & 0 & 0 \\ 
        -1 & 2 & -1 & -1 \\ 
        -1 & 1 & 0 & -1 \\ 
        0 & 1 & 0 & -1 \\ 
     \end{bmatrix*}$
&     $\begin{array}{c}
        s_4s_3s_1s_2s_3 \\ 
     \end{array}$
   \\
\hline
     \end{tabular}
\caption{\footnotesize Weyl group $D_4$, level 5, elements 0-13}
\label{tab_level_5_1}
\end{table}
\begin{table}[H]
\centering
\renewcommand{\arraystretch}{0.98}
     \begin{tabular}{||c|c|c|c|c||}
     \hline
          &&&& \\ 
    $N^o$  & Weight &  Element   & Matrix & Inverse   \\ 
           &        & \& order  &        &           \\ 
     \hline
     \footnotesize
     $\begin{array}{c}
         14 \\ 
          \\ 
         (14) \\ 
     \end{array}$
&     \footnotesize
 $[-2, 1, -2, 4]$&     $\begin{array}{c}
        s_3s_2s_1s_2s_3 \\ 
             \end{array}$
&     \footnotesize
     $\begin{bmatrix*}[r]
        0 & 0 & -1 & 1 \\ 
        -1 & 1 & -1 & 1 \\ 
        -1 & 0 & 0 & 1 \\ 
        0 & 0 & 0 & 1 \\ 
     \end{bmatrix*}$
&     $\begin{array}{c}
        s_3s_2s_1s_2s_3 \\ 
     \end{array}$
   \\ 
     \hline
     \footnotesize
     $\begin{array}{c}
         15 \\ 
          \\ 
         (23) \\ 
     \end{array}$
&     \footnotesize
 $[-2, 3, 2, -4]$&     $\begin{array}{c}
        s_4s_2s_1s_2s_3 \\ 
             \end{array}$
&     \footnotesize
     $\begin{bmatrix*}[r]
        0 & 0 & -1 & 1 \\ 
        -1 & 1 & -1 & 1 \\ 
        0 & 1 & -1 & 0 \\ 
        -1 & 1 & -1 & 0 \\ 
     \end{bmatrix*}$
&     $\begin{array}{c}
        s_3s_2s_1s_2s_4 \\ 
     \end{array}$
   \\ 
     \hline
     \footnotesize
     $\begin{array}{c}
         16 \\ 
          \\ 
         (6) \\ 
     \end{array}$
&     \footnotesize
 $[-1, -2, 1, 5]$&     $\begin{array}{c}
        s_2s_3s_1s_2s_3 \\ 
             \end{array}$
&     \footnotesize
     $\begin{bmatrix*}[r]
        0 & 0 & -1 & 1 \\ 
        0 & -1 & 0 & 2 \\ 
        1 & -1 & 0 & 1 \\ 
        0 & 0 & 0 & 1 \\ 
     \end{bmatrix*}$
&     $\begin{array}{c}
        s_3s_2s_3s_1s_2 \\ 
     \end{array}$
   \\
\hline
     \end{tabular}
\caption{\footnotesize Weyl group $D_4$, level 5, elements 14-27}
\label{tab_level_5_2}
\end{table}

\begin{table}[H]
\centering
\renewcommand{\arraystretch}{0.98}
     \begin{tabular}{||c|c|c|c|c||}
     \hline
          &&&& \\ 
    $N^o$  & Weight &  Element   & Matrix & Inverse   \\ 
           &        & \& order  &        &           \\ 
     \hline
     \footnotesize
     $\begin{array}{c}
         0 \\ 
          \\ 
         (0) \\ 
     \end{array}$
&     \footnotesize
 $[-5, 1, 1, 1]$&     $\begin{array}{c}
        s_1s_2s_4s_3s_2s_1 \\ 
             \end{array}$
&     \footnotesize
     $\begin{bmatrix*}[r]
        -1 & 0 & 0 & 0 \\ 
        -2 & 1 & 0 & 0 \\ 
        -1 & 0 & 0 & 1 \\ 
        -1 & 0 & 1 & 0 \\ 
     \end{bmatrix*}$
&     $\begin{array}{c}
        s_1s_2s_4s_3s_2s_1 \\ 
     \end{array}$
   \\ 
     \hline
     \footnotesize
     $\begin{array}{c}
         1 \\ 
          \\ 
         (11) \\ 
     \end{array}$
&     \footnotesize
 $[5, -3, -1, 1]$&     $\begin{array}{c}
        s_3s_2s_4s_3s_2s_1 \\ 
             \end{array}$
&     \footnotesize
     $\begin{bmatrix*}[r]
        -1 & 1 & 0 & 0 \\ 
        -2 & 1 & 0 & 0 \\ 
        -1 & 1 & 0 & -1 \\ 
        -1 & 0 & 1 & 0 \\ 
     \end{bmatrix*}$
&     $\begin{array}{c}
        s_4s_1s_2s_4s_3s_2 \\ 
     \end{array}$
   \\ 
     \hline
     \footnotesize
     $\begin{array}{c}
         2 \\ 
          \\ 
         (10) \\ 
     \end{array}$
&     \footnotesize
 $[5, -3, 1, -1]$&     $\begin{array}{c}
        s_4s_2s_4s_3s_2s_1 \\ 
             \end{array}$
&     \footnotesize
     $\begin{bmatrix*}[r]
        -1 & 1 & 0 & 0 \\ 
        -2 & 1 & 0 & 0 \\ 
        -1 & 0 & 0 & 1 \\ 
        -1 & 1 & -1 & 0 \\ 
     \end{bmatrix*}$
&     $\begin{array}{c}
        s_3s_1s_2s_4s_3s_2 \\ 
     \end{array}$
   \\
\hline
     \end{tabular}
\caption{\footnotesize Weyl group $D_4$, level 6, elements 0-13}
\label{tab_level_6_1}
\end{table}
\begin{table}[H]
\centering
\renewcommand{\arraystretch}{0.98}
     \begin{tabular}{||c|c|c|c|c||}
     \hline
          &&&& \\ 
    $N^o$  & Weight &  Element   & Matrix & Inverse   \\ 
           &        & \& order  &        &           \\ 
     \hline
     \footnotesize
     $\begin{array}{c}
         14 \\ 
          \\ 
         (14) \\ 
     \end{array}$
&     \footnotesize
 $[-1, -1, -1, 5]$&     $\begin{array}{c}
        s_3s_2s_3s_1s_2s_3 \\ 
             \end{array}$
&     \footnotesize
     $\begin{bmatrix*}[r]
        0 & 0 & -1 & 1 \\ 
        0 & -1 & 0 & 2 \\ 
        -1 & 0 & 0 & 1 \\ 
        0 & 0 & 0 & 1 \\ 
     \end{bmatrix*}$
&     $\begin{array}{c}
        s_3s_2s_3s_1s_2s_3 \\ 
     \end{array}$
   \\ 
     \hline
     \footnotesize
     $\begin{array}{c}
         15 \\ 
          \\ 
         (25) \\ 
     \end{array}$
&     \footnotesize
 $[-1, 3, 1, -5]$&     $\begin{array}{c}
        s_4s_2s_3s_1s_2s_3 \\ 
             \end{array}$
&     \footnotesize
     $\begin{bmatrix*}[r]
        0 & 0 & -1 & 1 \\ 
        0 & -1 & 0 & 2 \\ 
        1 & -1 & 0 & 1 \\ 
        0 & -1 & 0 & 1 \\ 
     \end{bmatrix*}$
&     $\begin{array}{c}
        s_3s_2s_3s_1s_2s_4 \\ 
     \end{array}$
   \\ 
     \hline
     \footnotesize
     $\begin{array}{c}
         16 \\ 
          \\ 
         (6) \\ 
     \end{array}$
&     \footnotesize
 $[2, -5, 4, 2]$&     $\begin{array}{c}
        s_2s_4s_3s_1s_2s_3 \\ 
             \end{array}$
&     \footnotesize
     $\begin{bmatrix*}[r]
        0 & 0 & -1 & 1 \\ 
        1 & -1 & -1 & 1 \\ 
        1 & -1 & 0 & 1 \\ 
        1 & 0 & -1 & 0 \\ 
     \end{bmatrix*}$
&     $\begin{array}{c}
        s_3s_2s_4s_3s_1s_2 \\ 
     \end{array}$
   \\
\hline
     \end{tabular}
\caption{\footnotesize Weyl group $D_4$, level 6, elements 14-27}
\label{tab_level_6_2}
\end{table}
\begin{table}[H]
\centering
\renewcommand{\arraystretch}{0.98}
     \begin{tabular}{||c|c|c|c|c||}
     \hline
          &&&& \\ 
    $N^o$  & Weight &  Element   & Matrix & Inverse   \\ 
           &        & \& order  &        &           \\ 
     \hline
     \footnotesize
     $\begin{array}{c}
         28 \\ 
          \\ 
         (18) \\ 
     \end{array}$
&     \footnotesize
 $[-1, 3, -5, 1]$&     $\begin{array}{c}
        s_3s_2s_4s_1s_2s_4 \\ 
             \end{array}$
&     \footnotesize
     $\begin{bmatrix*}[r]
        0 & 0 & 1 & -1 \\ 
        0 & -1 & 2 & 0 \\ 
        0 & -1 & 1 & 0 \\ 
        1 & -1 & 1 & 0 \\ 
     \end{bmatrix*}$
&     $\begin{array}{c}
        s_4s_2s_4s_1s_2s_3 \\ 
     \end{array}$
   \\ 
     \hline
     \footnotesize
     $\begin{array}{c}
         29 \\ 
          \\ 
         (29) \\ 
     \end{array}$
&     \footnotesize
 $[-1, -1, 5, -1]$&     $\begin{array}{c}
        s_4s_2s_4s_1s_2s_4 \\ 
             \end{array}$
&     \footnotesize
     $\begin{bmatrix*}[r]
        0 & 0 & 1 & -1 \\ 
        0 & -1 & 2 & 0 \\ 
        0 & 0 & 1 & 0 \\ 
        -1 & 0 & 1 & 0 \\ 
     \end{bmatrix*}$
&     $\begin{array}{c}
        s_4s_2s_4s_1s_2s_4 \\ 
     \end{array}$
   \\ 
     \hline
     \end{tabular}
\caption{\footnotesize Weyl group $D_4$, level 6, elements 28-29}
\label{tab_level_6_}
\end{table}

\begin{table}[H]
\centering
\renewcommand{\arraystretch}{0.98}
     \begin{tabular}{||c|c|c|c|c||}
     \hline
          &&&& \\ 
    $N^o$  & Weight &  Element   & Matrix & Inverse   \\ 
           &        & \& order  &        &           \\ 
     \hline
     \footnotesize
     $\begin{array}{c}
         0 \\ 
          \\ 
         (0) \\ 
     \end{array}$
&     \footnotesize
 $[-4, -1, 2, 2]$&     $\begin{array}{c}
        s_2s_1s_2s_4s_3s_2s_1 \\ 
             \end{array}$
&     \footnotesize
     $\begin{bmatrix*}[r]
        -1 & 0 & 0 & 0 \\ 
        -1 & -1 & 1 & 1 \\ 
        -1 & 0 & 0 & 1 \\ 
        -1 & 0 & 1 & 0 \\ 
     \end{bmatrix*}$
&     $\begin{array}{c}
        s_2s_1s_2s_4s_3s_2s_1 \\ 
     \end{array}$
   \\ 
     \hline
     \footnotesize
     $\begin{array}{c}
         1 \\ 
          \\ 
         (2) \\ 
     \end{array}$
&     \footnotesize
 $[-5, 2, -1, 1]$&     $\begin{array}{c}
        s_3s_1s_2s_4s_3s_2s_1 \\ 
             \end{array}$
&     \footnotesize
     $\begin{bmatrix*}[r]
        -1 & 0 & 0 & 0 \\ 
        -2 & 1 & 0 & 0 \\ 
        -1 & 1 & 0 & -1 \\ 
        -1 & 0 & 1 & 0 \\ 
     \end{bmatrix*}$
&     $\begin{array}{c}
        s_4s_1s_2s_4s_3s_2s_1 \\ 
     \end{array}$
   \\ 
     \hline
     \footnotesize
     $\begin{array}{c}
         2 \\ 
          \\ 
         (1) \\ 
     \end{array}$
&     \footnotesize
 $[-5, 2, 1, -1]$&     $\begin{array}{c}
        s_4s_1s_2s_4s_3s_2s_1 \\ 
             \end{array}$
&     \footnotesize
     $\begin{bmatrix*}[r]
        -1 & 0 & 0 & 0 \\ 
        -2 & 1 & 0 & 0 \\ 
        -1 & 0 & 0 & 1 \\ 
        -1 & 1 & -1 & 0 \\ 
     \end{bmatrix*}$
&     $\begin{array}{c}
        s_3s_1s_2s_4s_3s_2s_1 \\ 
     \end{array}$
   \\
\hline
     \end{tabular}
\caption{\footnotesize Weyl group $D_4$, level 7, elements 0-13}
\label{tab_level_7_1}
\end{table}
\begin{table}[H]
\centering
\renewcommand{\arraystretch}{0.98}
     \begin{tabular}{||c|c|c|c|c||}
     \hline
          &&&& \\ 
    $N^o$  & Weight &  Element   & Matrix & Inverse   \\ 
           &        & \& order  &        &           \\ 
     \hline
     \footnotesize
     $\begin{array}{c}
         14 \\ 
          \\ 
         (14) \\ 
     \end{array}$
&     \footnotesize
 $[2, -1, -4, 2]$&     $\begin{array}{c}
        s_3s_2s_4s_3s_1s_2s_3 \\ 
             \end{array}$
&     \footnotesize
     $\begin{bmatrix*}[r]
        0 & 0 & -1 & 1 \\ 
        1 & -1 & -1 & 1 \\ 
        0 & 0 & -1 & 0 \\ 
        1 & 0 & -1 & 0 \\ 
     \end{bmatrix*}$
&     $\begin{array}{c}
        s_3s_2s_4s_3s_1s_2s_3 \\ 
     \end{array}$
   \\ 
     \hline
     \footnotesize
     $\begin{array}{c}
         15 \\ 
          \\ 
         (25) \\ 
     \end{array}$
&     \footnotesize
 $[2, -3, 4, -2]$&     $\begin{array}{c}
        s_4s_2s_4s_3s_1s_2s_3 \\ 
             \end{array}$
&     \footnotesize
     $\begin{bmatrix*}[r]
        0 & 0 & -1 & 1 \\ 
        1 & -1 & -1 & 1 \\ 
        1 & -1 & 0 & 1 \\ 
        0 & -1 & 0 & 1 \\ 
     \end{bmatrix*}$
&     $\begin{array}{c}
        s_3s_2s_4s_3s_1s_2s_4 \\ 
     \end{array}$
   \\ 
     \hline
     \footnotesize
     $\begin{array}{c}
         16 \\ 
          \\ 
         (21) \\ 
     \end{array}$
&     \footnotesize
 $[1, 2, -5, -1]$&     $\begin{array}{c}
        s_4s_3s_2s_4s_1s_2s_3 \\ 
             \end{array}$
&     \footnotesize
     $\begin{bmatrix*}[r]
        0 & 0 & -1 & 1 \\ 
        0 & 1 & -2 & 0 \\ 
        0 & 0 & -1 & 0 \\ 
        -1 & 1 & -1 & 0 \\ 
     \end{bmatrix*}$
&     $\begin{array}{c}
        s_3s_2s_4s_1s_2s_4s_3 \\ 
     \end{array}$
   \\
\hline
     \end{tabular}
\caption{\footnotesize Weyl group $D_4$, level 7, elements 14-27}
\label{tab_level_7_2}
\end{table}

\begin{table}[H]
\centering
\renewcommand{\arraystretch}{0.98}
     \begin{tabular}{||c|c|c|c|c||}
     \hline
          &&&& \\ 
    $N^o$  & Weight &  Element   & Matrix & Inverse   \\ 
           &        & \& order  &        &           \\ 
     \hline
     \footnotesize
     $\begin{array}{c}
         0 \\ 
          \\ 
         (4) \\ 
     \end{array}$
&     \footnotesize
 $[-4, 1, -2, 2]$&     $\begin{array}{c}
        s_3s_2s_1s_2s_4s_3s_2s_1 \\ 
             \end{array}$
&     \footnotesize
     $\begin{bmatrix*}[r]
        -1 & 0 & 0 & 0 \\ 
        -1 & -1 & 1 & 1 \\ 
        0 & -1 & 1 & 0 \\ 
        -1 & 0 & 1 & 0 \\ 
     \end{bmatrix*}$
&     $\begin{array}{c}
        s_2s_4s_1s_2s_4s_3s_2s_1 \\ 
     \end{array}$
   \\ 
     \hline
     \footnotesize
     $\begin{array}{c}
         1 \\ 
          \\ 
         (2) \\ 
     \end{array}$
&     \footnotesize
 $[-4, 1, 2, -2]$&     $\begin{array}{c}
        s_4s_2s_1s_2s_4s_3s_2s_1 \\ 
             \end{array}$
&     \footnotesize
     $\begin{bmatrix*}[r]
        -1 & 0 & 0 & 0 \\ 
        -1 & -1 & 1 & 1 \\ 
        -1 & 0 & 0 & 1 \\ 
        0 & -1 & 0 & 1 \\ 
     \end{bmatrix*}$
&     $\begin{array}{c}
        s_2s_3s_1s_2s_4s_3s_2s_1 \\ 
     \end{array}$
   \\ 
     \hline
     \footnotesize
     $\begin{array}{c}
         2 \\ 
          \\ 
         (1) \\ 
     \end{array}$
&     \footnotesize
 $[-3, -2, 1, 3]$&     $\begin{array}{c}
        s_2s_3s_1s_2s_4s_3s_2s_1 \\ 
             \end{array}$
&     \footnotesize
     $\begin{bmatrix*}[r]
        -1 & 0 & 0 & 0 \\ 
        -1 & 0 & 1 & -1 \\ 
        -1 & 1 & 0 & -1 \\ 
        -1 & 0 & 1 & 0 \\ 
     \end{bmatrix*}$
&     $\begin{array}{c}
        s_4s_2s_1s_2s_4s_3s_2s_1 \\ 
     \end{array}$
   \\
\hline
     \end{tabular}
\caption{\footnotesize Weyl group $D_4$, level 8, elements 0-13}
\label{tab_level_8_1}
\end{table}
\begin{table}[H]
\centering
\renewcommand{\arraystretch}{0.98}
     \begin{tabular}{||c|c|c|c|c||}
     \hline
          &&&& \\ 
    $N^o$  & Weight &  Element   & Matrix & Inverse   \\ 
           &        & \& order  &        &           \\ 
     \hline
     \footnotesize
     $\begin{array}{c}
         14 \\ 
          \\ 
         (17) \\ 
     \end{array}$
&     \footnotesize
 $[2, 1, -4, -2]$&     $\begin{array}{c}
        s_4s_3s_2s_4s_3s_1s_2s_3 \\ 
             \end{array}$
&     \footnotesize
     $\begin{bmatrix*}[r]
        0 & 0 & -1 & 1 \\ 
        1 & -1 & -1 & 1 \\ 
        0 & 0 & -1 & 0 \\ 
        0 & -1 & 0 & 1 \\ 
     \end{bmatrix*}$
&     $\begin{array}{c}
        s_3s_2s_4s_3s_1s_2s_4s_3 \\ 
     \end{array}$
   \\ 
     \hline
     \footnotesize
     $\begin{array}{c}
         15 \\ 
          \\ 
         (6) \\ 
     \end{array}$
&     \footnotesize
 $[2, -5, 2, 2]$&     $\begin{array}{c}
        s_2s_4s_3s_2s_1s_2s_4s_3 \\ 
             \end{array}$
&     \footnotesize
     $\begin{bmatrix*}[r]
        0 & 1 & -1 & -1 \\ 
        -1 & 1 & -1 & -1 \\ 
        -1 & 1 & 0 & -1 \\ 
        -1 & 1 & -1 & 0 \\ 
     \end{bmatrix*}$
&     $\begin{array}{c}
        s_4s_3s_2s_1s_2s_4s_3s_2 \\ 
     \end{array}$
   \\ 
     \hline
     \footnotesize
     $\begin{array}{c}
         16 \\ 
          \\ 
         (16) \\ 
     \end{array}$
&     \footnotesize
 $[-1, 3, -1, -5]$&     $\begin{array}{c}
        s_4s_3s_2s_3s_1s_2s_4s_3 \\ 
             \end{array}$
&     \footnotesize
     $\begin{bmatrix*}[r]
        0 & 1 & -1 & -1 \\ 
        0 & 1 & 0 & -2 \\ 
        -1 & 1 & 0 & -1 \\ 
        0 & 0 & 0 & -1 \\ 
     \end{bmatrix*}$
&     $\begin{array}{c}
        s_4s_3s_2s_3s_1s_2s_4s_3 \\ 
     \end{array}$
   \\
\hline
     \end{tabular}
\caption{\footnotesize Weyl group $D_4$, level 8, elements 14-22}
\label{tab_level_8_2}
\end{table}

\begin{table}[H]
\centering
\renewcommand{\arraystretch}{0.98}
     \begin{tabular}{||c|c|c|c|c||}
     \hline
          &&&& \\ 
    $N^o$  & Weight &  Element   & Matrix & Inverse   \\ 
           &        & \& order  &        &           \\ 
     \hline
     \footnotesize
     $\begin{array}{c}
         0 \\ 
          \\ 
         (3) \\ 
     \end{array}$
&     \footnotesize
 $[-4, 3, -2, -2]$&     $\begin{array}{c}
        s_4s_3s_2s_1s_2s_4s_3s_2s_1 \\ 
             \end{array}$
&     \footnotesize
     $\begin{bmatrix*}[r]
        -1 & 0 & 0 & 0 \\ 
        -1 & -1 & 1 & 1 \\ 
        0 & -1 & 1 & 0 \\ 
        0 & -1 & 0 & 1 \\ 
     \end{bmatrix*}$
&     $\begin{array}{c}
        s_2s_4s_3s_1s_2s_4s_3s_2s_1 \\ 
     \end{array}$
   \\ 
     \hline
     \footnotesize
     $\begin{array}{c}
         1 \\ 
          \\ 
         (5) \\ 
     \end{array}$
&     \footnotesize
 $[-3, -1, -1, 3]$&     $\begin{array}{c}
        s_3s_2s_3s_1s_2s_4s_3s_2s_1 \\ 
             \end{array}$
&     \footnotesize
     $\begin{bmatrix*}[r]
        -1 & 0 & 0 & 0 \\ 
        -1 & 0 & 1 & -1 \\ 
        0 & -1 & 1 & 0 \\ 
        -1 & 0 & 1 & 0 \\ 
     \end{bmatrix*}$
&     $\begin{array}{c}
        s_4s_2s_4s_1s_2s_4s_3s_2s_1 \\ 
     \end{array}$
   \\ 
     \hline
     \footnotesize
     $\begin{array}{c}
         2 \\ 
          \\ 
         (2) \\ 
     \end{array}$
&     \footnotesize
 $[-3, 1, 1, -3]$&     $\begin{array}{c}
        s_4s_2s_3s_1s_2s_4s_3s_2s_1 \\ 
             \end{array}$
&     \footnotesize
     $\begin{bmatrix*}[r]
        -1 & 0 & 0 & 0 \\ 
        -1 & 0 & 1 & -1 \\ 
        -1 & 1 & 0 & -1 \\ 
        0 & 0 & 0 & -1 \\ 
     \end{bmatrix*}$
&     $\begin{array}{c}
        s_4s_2s_3s_1s_2s_4s_3s_2s_1 \\ 
     \end{array}$
   \\
\hline
     \end{tabular}
\caption{\footnotesize Weyl group $D_4$, level 9, elements 0-13}
\label{tab_level_9_1}
\end{table}
\begin{table}[H]
\centering
\renewcommand{\arraystretch}{0.98}
     \begin{tabular}{||c|c|c|c|c||}
     \hline
          &&&& \\ 
    $N^o$  & Weight &  Element   & Matrix & Inverse   \\ 
           &        & \& order  &        &           \\ 
     \hline
     \footnotesize
     $\begin{array}{c}
         14 \\ 
          \\ 
         (14) \\ 
     \end{array}$
&     \footnotesize
 $[1, 1, -3, -3]$&     $\begin{array}{c}
        s_4s_3s_2s_4s_3s_1s_2s_4s_3 \\ 
             \end{array}$
&     \footnotesize
     $\begin{bmatrix*}[r]
        0 & 1 & -1 & -1 \\ 
        1 & 0 & -1 & -1 \\ 
        0 & 0 & -1 & 0 \\ 
        0 & 0 & 0 & -1 \\ 
     \end{bmatrix*}$
&     $\begin{array}{c}
        s_4s_3s_2s_4s_3s_1s_2s_4s_3 \\ 
     \end{array}$
   \\ 
     \hline
     \footnotesize
     $\begin{array}{c}
         15 \\ 
          \\ 
         (9) \\ 
     \end{array}$
&     \footnotesize
 $[3, -1, -1, -3]$&     $\begin{array}{c}
        s_4s_3s_2s_4s_3s_2s_1s_2s_4 \\ 
             \end{array}$
&     \footnotesize
     $\begin{bmatrix*}[r]
        0 & 0 & 1 & -1 \\ 
        -1 & 0 & 1 & -1 \\ 
        0 & -1 & 1 & 0 \\ 
        0 & 0 & 0 & -1 \\ 
     \end{bmatrix*}$
&     $\begin{array}{c}
        s_4s_2s_4s_3s_1s_2s_4s_3s_2 \\ 
     \end{array}$
   \\ 
     \hline
     \end{tabular}
\caption{\footnotesize Weyl group $D_4$, level 9, elements 14-15}
\label{tab_level_9_2}
\end{table}

\begin{table}[H]
\centering
\renewcommand{\arraystretch}{0.98}
     \begin{tabular}{||c|c|c|c|c||}
     \hline
          &&&& \\ 
    $N^o$  & Weight &  Element   & Matrix & Inverse   \\ 
           &        & \& order  &        &           \\ 
     \hline
     \footnotesize
     $\begin{array}{c}
         0 \\ 
          \\ 
         (0) \\ 
     \end{array}$
&     \footnotesize
 $[-1, -3, 1, 1]$&     $\begin{array}{c}
        s_2s_4s_3s_2s_1s_2s_4s_3s_2s_1 \\ 
             \end{array}$
&     \footnotesize
     $\begin{bmatrix*}[r]
        -1 & 0 & 0 & 0 \\ 
        0 & -1 & 0 & 0 \\ 
        0 & -1 & 1 & 0 \\ 
        0 & -1 & 0 & 1 \\ 
     \end{bmatrix*}$
&     $\begin{array}{c}
        s_2s_4s_3s_2s_1s_2s_4s_3s_2s_1 \\ 
     \end{array}$
   \\ 
     \hline
     \footnotesize
     $\begin{array}{c}
         1 \\ 
          \\ 
         (3) \\ 
     \end{array}$
&     \footnotesize
 $[-3, 2, -1, -3]$&     $\begin{array}{c}
        s_4s_3s_2s_3s_1s_2s_4s_3s_2s_1 \\ 
             \end{array}$
&     \footnotesize
     $\begin{bmatrix*}[r]
        -1 & 0 & 0 & 0 \\ 
        -1 & 0 & 1 & -1 \\ 
        0 & -1 & 1 & 0 \\ 
        0 & 0 & 0 & -1 \\ 
     \end{bmatrix*}$
&     $\begin{array}{c}
        s_4s_2s_4s_3s_1s_2s_4s_3s_2s_1 \\ 
     \end{array}$
   \\ 
     \hline
     \footnotesize
     $\begin{array}{c}
         2 \\ 
          \\ 
         (4) \\ 
     \end{array}$
&     \footnotesize
 $[-2, -1, -2, 2]$&     $\begin{array}{c}
        s_3s_2s_4s_3s_1s_2s_4s_3s_2s_1 \\ 
             \end{array}$
&     \footnotesize
     $\begin{bmatrix*}[r]
        -1 & 0 & 0 & 0 \\ 
        -1 & 1 & -1 & -1 \\ 
        0 & 0 & -1 & 0 \\ 
        -1 & 1 & -1 & 0 \\ 
     \end{bmatrix*}$
&     $\begin{array}{c}
        s_4s_3s_2s_4s_1s_2s_4s_3s_2s_1 \\ 
     \end{array}$
   \\
\hline
     \end{tabular}
\caption{\footnotesize Weyl group $D_4$, level 10, elements 0-8}
\label{tab_level_10_1}
\end{table}

\begin{table}[H]
\centering
\renewcommand{\arraystretch}{0.98}
     \begin{tabular}{||c|c|c|c|c||}
     \hline
    $N^o$  & Weight &  Element   & Matrix & Inverse   \\
           &        & \& order  &        &           \\
     \hline
     \footnotesize
     $\begin{array}{c}
         0 \\
          \\
         (0) \\
     \end{array}$
&     \footnotesize
 $-1,-2,-1,1$&     $\begin{array}{c}
        s_3s_2s_4s_3s_2s_1s_2s_4s_3s_2s_1 \\
                2 \\
     \end{array}$
&     \footnotesize
     $\begin{bmatrix*}[r]
        -1 & 0 & 0 & 0 \\
        0 & -1 & 0 & 0 \\
        0 & 0 & -1 & 0 \\
        0 & -1 & 0 & 1 \\
     \end{bmatrix*}$
&     $\begin{array}{c}
        s_3s_2s_4s_3s_2s_1s_2s_4s_3s_2s_1 \\
     \end{array}$
   \\
     \hline
     \footnotesize
     $\begin{array}{c}
         1 \\
          \\
         (1) \\
     \end{array}$
&     \footnotesize
 $-1,-2,1,-1$&     $\begin{array}{c}
        s_4s_2s_4s_3s_2s_1s_2s_4s_3s_2s_1 \\
                2 \\
     \end{array}$
&     \footnotesize
     $\begin{bmatrix*}[r]
        -1 & 0 & 0 & 0 \\
        0 & -1 & 0 & 0 \\
        0 & -1 & 1 & 0 \\
        0 & 0 & 0 & -1 \\
     \end{bmatrix*}$
&     $\begin{array}{c}
        s_4s_2s_4s_3s_2s_1s_2s_4s_3s_2s_1 \\
     \end{array}$
   \\
     \hline
     \footnotesize
     $\begin{array}{c}
         2 \\
          \\
         (2) \\
     \end{array}$
&     \footnotesize
 $-2,1,-2,-2$&     $\begin{array}{c}
        s_4s_3s_2s_4s_3s_1s_2s_4s_3s_2s_1 \\
                2 \\
     \end{array}$
&     \footnotesize
     $\begin{bmatrix*}[r]
        -1 & 0 & 0 & 0 \\
        -1 & 1 & -1 & -1 \\
        0 & 0 & -1 & 0 \\
        0 & 0 & 0 & -1 \\
     \end{bmatrix*}$
&     $\begin{array}{c}
        s_4s_3s_2s_4s_3s_1s_2s_4s_3s_2s_1 \\
     \end{array}$
   \\
\hline
     \end{tabular}
\caption{\footnotesize Weyl group $D_4$, level 11, elements 0-3}
\label{tab_level_11_1}

\begin{tabular}{||c|c|c|c|c||}
     \hline
    $N^o$  & Weight &  Element   & Matrix & Inverse   \\
           &        & \& order  &        &           \\
     \hline
     \footnotesize
     $\begin{array}{c}
         0 \\
          \\
         (0) \\
     \end{array}$
&     \footnotesize
 $-1,-1,-1,-1$&     $\begin{array}{c}
        s_4s_3s_2s_4s_3s_2s_1s_2s_4s_3s_2s_1 \\
                2 \\
     \end{array}$
&     \footnotesize
     $\begin{bmatrix*}[r]
        -1 & 0 & 0 & 0 \\
        0 & -1 & 0 & 0 \\
        0 & 0 & -1 & 0 \\
        0 & 0 & 0 & -1 \\
     \end{bmatrix*}$
&     $\begin{array}{c}
      \text{itself}
     \end{array}$
   \\
     \hline
\end{tabular}

\caption{\footnotesize Weyl group $D_4$, level 12, one element}
\label{tab_level_11_12}
\end{table}

~\\
\section{\bf Thirteen conjugacy classes of $W(D_4)$}
 \label{sec_conj_classes}

\begin{table}[H]
\centering
\renewcommand{\arraystretch}{0.98}
\begin{tabular}{||c|c|c|c||}
\hline
 $N^o$ in CCL  &  Element & Level &  $N^o$ in Level    \\
\hline
\footnotesize 1 & $e$ & \footnotesize 0 & \footnotesize 1 \\
\hline
\end{tabular}
\caption{\footnotesize Weyl group $D_4$, conjugacy class 0,  order =  1}

\begin{tabular}{||c|c|c|c||}
\hline
  $N^o$ in CCL  & Element & Level & $N^o$ in Level    \\
\hline
\footnotesize 1 & $s_1$ & \footnotesize 1 & \footnotesize 1 \\
\footnotesize 2 & $s_2$ & \footnotesize 1 & \footnotesize 2 \\
\footnotesize 3 & $s_3$ & \footnotesize 1 & \footnotesize 3 \\
\footnotesize 4 & $s_4$ & \footnotesize 1 & \footnotesize 4 \\
\footnotesize 5 & $s_2s_1s_2$ & \footnotesize 3 & \footnotesize 6 \\
\footnotesize 6 & $s_3s_2s_3$ & \footnotesize 3 & \footnotesize 11 \\
\footnotesize 7 & $s_4s_2s_4$ & \footnotesize 3 & \footnotesize 16 \\
\footnotesize 8 & $s_3s_2s_1s_2s_3$ & \footnotesize 5 & \footnotesize 15 \\
\footnotesize 9 & $s_4s_3s_2s_4s_3$ & \footnotesize 5 & \footnotesize 23 \\
\footnotesize 10 & $s_4s_2s_1s_2s_4$ & \footnotesize 5 & \footnotesize 25 \\
\footnotesize 11 & $s_4s_3s_2s_1s_2s_4s_3$ & \footnotesize 7 & \footnotesize 18 \\
\footnotesize 12 & $s_2s_4s_3s_2s_1s_2s_4s_3s_2$ & \footnotesize 9 & \footnotesize 7 \\
\hline
\end{tabular}
\caption{\footnotesize Weyl group $D_4$, conjugacy class 1,  order =  2}

\label{tableCCL_0123}
\end{table} 

\begin{table}[H]
\centering
\renewcommand{\arraystretch}{0.98}
\begin{tabular}{||c|c|c|c||}
\hline
          &&& \\
 $N^o$ in CCL  & Element & Level & $N^o$ in Level    \\
          &&& \\
\hline
\footnotesize 1 & $s_2s_1$ & \footnotesize 2 & \footnotesize 0 \\
\footnotesize 2 & $s_1s_2$ & \footnotesize 2 & \footnotesize 3 \\
\footnotesize 3 & $s_3s_2$ & \footnotesize 2 & \footnotesize 4 \\
\footnotesize 4 & $s_4s_2$ & \footnotesize 2 & \footnotesize 5 \\
\footnotesize 5 & $s_2s_3$ & \footnotesize 2 & \footnotesize 6 \\
\hline
\end{tabular}
\vspace{2mm}\caption{\footnotesize Weyl group $D_4$, conjugacy class 2,  order =  3}
\label{tableCCL_2}
\end{table} 

\begin{table}[H]
\centering
\renewcommand{\arraystretch}{0.98}

\begin{tabular}{||c|c|c|c||}
\hline
          &&& \\
 $N^o$ in CCL  & Element & Level & $N^o$ in Level    \\
          &&& \\
\hline
\footnotesize 1 & $s_3s_1$ & \footnotesize 2 & \footnotesize 1 \\
\footnotesize 2 & $s_2s_3s_1s_2$ & \footnotesize 4 & \footnotesize 8 \\
\footnotesize 3 & $s_3s_2s_3s_1s_2s_3$ & \footnotesize 6 & \footnotesize 14 \\
\footnotesize 4 & $s_4s_2s_3s_1s_2s_4$ & \footnotesize 6 & \footnotesize 26 \\
\footnotesize 5 & $s_4s_3s_2s_3s_1s_2s_4s_3$ & \footnotesize 8 & \footnotesize 16 \\
\footnotesize 6 & $s_4s_2s_4s_3s_2s_1s_2s_4s_3s_2$ & \footnotesize 10 & \footnotesize 6 \\
\hline
\end{tabular}
\label{tableCCL_3}
\caption{\footnotesize Weyl group $D_4$, conjugacy class 3,  order =  2}
~\\
~\\
~\\
~\\
\begin{tabular}{||c|c|c|c||}
\hline
          &&& \\
 $N^o$ in CCL  & Element & Level & $N^o$ in Level    \\
          &&& \\
\hline
\footnotesize 1 & $s_4s_1$ & \footnotesize 2 & \footnotesize 2 \\
\footnotesize 2 & $s_2s_4s_1s_2$ & \footnotesize 4 & \footnotesize 10 \\
\footnotesize 3 & $s_4s_2s_4s_1s_2s_4$ & \footnotesize 6 & \footnotesize 29 \\
\footnotesize 4 & $s_3s_2s_4s_1s_2s_3$ & \footnotesize 6 & \footnotesize 17 \\
\footnotesize 5 & $s_4s_3s_2s_4s_1s_2s_4s_3$ & \footnotesize 8 & \footnotesize 19 \\
\footnotesize 6 & $s_3s_2s_4s_3s_2s_1s_2s_4s_3s_2$ & \footnotesize 10 & \footnotesize 5 \\
\hline
\end{tabular}
\label{tableCCL_4}
\caption{\footnotesize Weyl group $D_4$, conjugacy class 4,  order =  2}
~\\
~\\
~\\
~\\
\begin{tabular}{||c|c|c|c||}
\hline
          &&& \\
 $N^o$ in CCL  & Element & Level & $N^o$ in Level    \\
          &&& \\
\hline
\footnotesize 1 & $s_4s_3$ & \footnotesize 2 & \footnotesize 7 \\
\footnotesize 2 & $s_2s_4s_3s_2$ & \footnotesize 4 & \footnotesize 11 \\
\footnotesize 3 & $s_1s_2s_4s_3s_2s_1$ & \footnotesize 6 & \footnotesize 0 \\
\footnotesize 4 & $s_4s_3s_2s_4s_3s_2$ & \footnotesize 6 & \footnotesize 12 \\
\footnotesize 5 & $s_4s_3s_1s_2s_4s_3s_2s_1$ & \footnotesize 8 & \footnotesize 3 \\
\footnotesize 6 & $s_2s_4s_3s_2s_1s_2s_4s_3s_2s_1$ & \footnotesize 10 & \footnotesize 0 \\
\hline
\end{tabular}
\vspace{2mm}\caption{\footnotesize Weyl group $D_4$, conjugacy class 5,  order =  2}
\label{tableCCL_345}
\end{table} 

\begin{table}[H]
\centering
\renewcommand{\arraystretch}{0.98}
\begin{tabular}{||c|c|c|c||}
\hline
          &&& \\
 $N^o$ in CCL  & Element & Level & $N^o$ in Level    \\
          &&& \\
\hline
\footnotesize 1 & $s_3s_2s_1$ & \footnotesize 3 & \footnotesize 0 \\
\footnotesize 2 & $s_2s_3s_1$ & \footnotesize 3 & \footnotesize 2 \\
\footnotesize 3 & $s_3s_1s_2$ & \footnotesize 3 & \footnotesize 6 \\
\footnotesize 4 & $s_1s_2s_3$ & \footnotesize 3 & \footnotesize 9 \\
\footnotesize 5 & $s_4s_2s_4s_3s_1$ & \footnotesize 5 & \footnotesize 3 \\
\hline
\end{tabular}
\vspace{2mm}\caption{\footnotesize Weyl group $D_4$, conjugacy class 6,  order =  4}
\label{tableCCL_6}
\end{table} 

\begin{table}[H]
\centering
\renewcommand{\arraystretch}{0.98}
\begin{tabular}{||c|c|c|c||}
\hline
          &&& \\
 $N^o$ in CCL  & Element & Level & $N^o$ in Level    \\
          &&& \\
\hline
\footnotesize 1 & $s_4s_2s_1$ & \footnotesize 3 & \footnotesize 1 \\
\footnotesize 2 & $s_2s_4s_1$ & \footnotesize 3 & \footnotesize 4 \\
\footnotesize 3 & $s_4s_1s_2$ & \footnotesize 3 & \footnotesize 7 \\
\footnotesize 4 & $s_1s_2s_4$ & \footnotesize 3 & \footnotesize 13 \\
\footnotesize 5 & $s_4s_3s_2s_3s_1$ & \footnotesize 5 & \footnotesize 1 \\
\hline
\end{tabular}
\vspace{2mm}\caption{\footnotesize Weyl group $D_4$, conjugacy class 7,  order =  4}
\label{tableCCL_7}
\end{table} 

\begin{table}[H]
\centering
\renewcommand{\arraystretch}{0.98}
\begin{tabular}{||c|c|c|c||}
\hline
          &&& \\
 $N^o$ in CCL  & Element & Level & $N^o$ in Level    \\
          &&& \\
\hline
\footnotesize 1 & $s_4s_3s_1$ & \footnotesize 3 & \footnotesize 3 \\
\footnotesize 2 & $s_2s_4s_3s_1s_2$ & \footnotesize 5 & \footnotesize 8 \\
\footnotesize 3 & $s_2s_1s_2s_4s_3s_2s_1$ & \footnotesize 7 & \footnotesize 0 \\
\footnotesize 4 & $s_3s_2s_4s_3s_1s_2s_3$ & \footnotesize 7 & \footnotesize 14 \\
\footnotesize 5 & $s_4s_2s_4s_3s_1s_2s_4$ & \footnotesize 7 & \footnotesize 26 \\
\footnotesize 6 & $s_4s_2s_3s_1s_2s_4s_3s_2s_1$ & \footnotesize 9 & \footnotesize 2 \\
\footnotesize 7 & $s_3s_2s_4s_1s_2s_4s_3s_2s_1$ & \footnotesize 9 & \footnotesize 4 \\
\footnotesize 8 & $s_4s_3s_2s_4s_3s_1s_2s_4s_3$ & \footnotesize 9 & \footnotesize 14 \\
\footnotesize 9 & $s_3s_2s_4s_3s_2s_1s_2s_4s_3s_2s_1$ & \footnotesize 11 & \footnotesize 0 \\
\footnotesize 10 & $s_4s_2s_4s_3s_2s_1s_2s_4s_3s_2s_1$ & \footnotesize 11 & \footnotesize 1 \\
\footnotesize 11 & $s_4s_3s_2s_4s_3s_1s_2s_4s_3s_2s_1$ & \footnotesize 11 & \footnotesize 2 \\
\footnotesize 12 & $s_4s_3s_2s_4s_3s_2s_1s_2s_4s_3s_2$ & \footnotesize 11 & \footnotesize 3 \\
\hline
\end{tabular}
\vspace{2mm}\caption{\footnotesize Weyl group $D_4$, conjugacy class 8,  order =  2}
\label{tableCCL_8}
\end{table} 

\begin{table}[H]
\centering
\renewcommand{\arraystretch}{0.98}
\begin{tabular}{||c|c|c|c||}
\hline
          &&& \\
 $N^o$ in CCL  & Element & Level & $N^o$ in Level    \\
          &&& \\
\hline
\footnotesize 1 & $s_4s_3s_2$ & \footnotesize 3 & \footnotesize 8 \\
\footnotesize 2 & $s_4s_2s_3$ & \footnotesize 3 & \footnotesize 11 \\
\footnotesize 3 & $s_2s_4s_3$ & \footnotesize 3 & \footnotesize 12 \\
\footnotesize 4 & $s_3s_2s_4$ & \footnotesize 3 & \footnotesize 14 \\
\footnotesize 5 & $s_2s_4s_3s_2s_1$ & \footnotesize 5 & \footnotesize 0 \\
\hline
\end{tabular}
\caption{\footnotesize Weyl group $D_4$, conjugacy class 9,  order =  4}
\label{tableCCL_9}
\end{table} 

\begin{table}[H]
\centering
\renewcommand{\arraystretch}{0.98}
\begin{tabular}{||c|c|c|c||}
\hline
          &&& \\
 $N^o$ in CCL  & Element & Level & $N^o$ in Level    \\
          &&& \\
\hline
\footnotesize 1 & $s_4s_3s_2s_1$ & \footnotesize 4 & \footnotesize 0 \\
\footnotesize 2 & $s_4s_2s_3s_1$ & \footnotesize 4 & \footnotesize 2 \\
\footnotesize 3 & $s_2s_4s_3s_1$ & \footnotesize 4 & \footnotesize 3 \\
\footnotesize 4 & $s_3s_2s_4s_1$ & \footnotesize 4 & \footnotesize 4 \\
\footnotesize 5 & $s_4s_3s_1s_2$ & \footnotesize 4 & \footnotesize 9 \\
\hline
\end{tabular}
\caption{\footnotesize Weyl group $D_4$, conjugacy class 10,  order =  6}
\label{tableCCL_10}
\end{table} 

\begin{table}[H]
\centering
\renewcommand{\arraystretch}{0.98}
\begin{tabular}{||c|c|c|c||}
\hline
          &&& \\
 $N^o$ in CCL  & Element & Level & $N^o$ in Level    \\
          &&& \\
\hline
\footnotesize 1 & $s_3s_2s_4s_3s_2s_1$ & \footnotesize 6 & \footnotesize 1 \\
\footnotesize 2 & $s_4s_2s_4s_3s_2s_1$ & \footnotesize 6 & \footnotesize 2 \\
\footnotesize 3 & $s_4s_3s_2s_3s_1s_2$ & \footnotesize 6 & \footnotesize 5 \\
\footnotesize 4 & $s_4s_3s_2s_4s_1s_2$ & \footnotesize 6 & \footnotesize 8 \\
\footnotesize 5 & $s_3s_1s_2s_4s_3s_2$ & \footnotesize 6 & \footnotesize 10 \\
\footnotesize 6 & $s_4s_1s_2s_4s_3s_2$ & \footnotesize 6 & \footnotesize 11 \\
\footnotesize 7 & $s_4s_2s_3s_1s_2s_3$ & \footnotesize 6 & \footnotesize 15 \\
\footnotesize 8 & $s_4s_2s_4s_1s_2s_3$ & \footnotesize 6 & \footnotesize 18 \\
\footnotesize 9 & $s_2s_3s_1s_2s_4s_3$ & \footnotesize 6 & \footnotesize 21 \\
\footnotesize 10 & $s_2s_4s_1s_2s_4s_3$ & \footnotesize 6 & \footnotesize 23 \\
\footnotesize 11 & $s_3s_2s_3s_1s_2s_4$ & \footnotesize 6 & \footnotesize 25 \\
\footnotesize 12 & $s_3s_2s_4s_1s_2s_4$ & \footnotesize 6 & \footnotesize 28 \\
\hline
\end{tabular}
\caption{\footnotesize Weyl group $D_4$, conjugacy class 11,  order =  4}
\label{tableCCL_11}

\begin{tabular}{||c|c|c|c||}
\hline
          &&& \\
 $N^o$ in CCL  & Element & Level & $N^o$ in Level    \\
          &&& \\
\hline
\footnotesize 1 & $s_4s_3s_2s_4s_3s_2s_1s_2s_4s_3s_2s_1$ & \footnotesize 12 & \footnotesize 0 \\
\hline
\end{tabular}
\caption{\footnotesize Weyl group $D_4$, conjugacy class 12,  order =  2}
\label{tableCCL_12}

\end{table} 

\begin{appendix}

\section{\bf Some properties of weights}
 \label{app_sect_weights}

This section lists some properties of finite Weyl groups,
weights related to Lie algebras and Weyl groups,
as well as some other concepts related to weights.

\subsection{Fundamental Weyl chamber}
   \label{sec_chambers}
Let $\varPhi$ be a root system, $\Delta$ be the set of the simple roots,
$\varPhi^+$ (resp. $\varPhi^-$) be the set of positive (resp. negative) roots,
$\Delta = \{\alpha_1,\dots,\alpha_l\}$,
and $\mathcal{E}$ be the linear space spanned by root of $\Delta$.
For any root $\alpha \in \varPhi$, let $H_{\alpha}$ be the hyperplane
$\{ x \in \mathcal{E} \mid (\alpha, x) = 0 \}$.
There are the finite number of the connected components of
\begin{equation*}
   \mathcal{E} - \bigcup\limits_{\alpha \in \varPhi} H_{\alpha}.
\end{equation*}
These components are called the open {\it Weyl chambers}.
There is the unique chamber $C$ such that for any $\xi \in C$, the following inequality holds:
\begin{equation}
   \label{eq_cond_fundam_domain}
   (\xi, \alpha_i) > 0 \text{ for all } \alpha_i \in \Delta,
\end{equation}
where $(\cdot,\cdot)$ is the Cartan-Killing bilinear form.
The unique Weyl chamber $C$ is called the {\it fundamental Weyl chamber}\footnotemark[1].

Eq. \eqref{eq_cond_fundam_domain} is equivalent to each of the following two statements:
\begin{equation}
   \label{eq_cond_fundam_domain_2}
   \begin{split}
   & (\xi, \alpha) > 0 \text{ for all } \alpha \in \varPhi^{+}, \\
   & (\xi, \alpha) < 0 \text{ for all } \alpha \in \varPhi^{-}. \\
   \end{split}
\end{equation}
\footnotetext[1]{The fundamental domains of usual and affine Weyl groups (Weyl chambers) were first described by E. Cartan in 1927, [C27], see \cite[p.62]{AR11}.}

\begin{theorem} {\cite[ch.VI, $\S$1, $n^\circ5$, Th.2]{Bo02}}
  \label{fact_1}
{\rm(i)} The Weyl group acts simply-transitively on the Weyl chambers.
Thus, the order of the Weyl group is equal to the number of Weyl chambers.

{\rm(ii)} Each $\xi \in \mathcal{E}$ is conjugate
to a unique point in the closure $\overline{C}$ of the fundamental Weyl chamber,
(i.e. $\overline{C}$  is a fundamental domain for $W$).
\end{theorem}
The word ``conjugate'' means ``in the same Weyl group orbit''.

\subsection{Dominant weights}
 \label{sec_dominant}
For any vectors $\alpha, \beta \in \varPhi$, let us define $\langle \alpha, \beta \rangle$
as follows:
\begin{equation}
  \label{eq_def_prod}
    \langle \alpha, \beta \rangle := \frac{2(\alpha, \beta)}{(\beta, \beta)}
\end{equation}
For the simply-laced Dynkin diagrams, if $\beta$ is a root then
$\langle \alpha, \beta \rangle  = (\alpha, \beta)$.
A {\it weight} (resp. {\it dominant weight}) is an element $\lambda \in \mathcal{E}$ such that
\begin{equation}
 \label{def_dominant}
    \langle \lambda, \alpha \rangle \in \mathbb{Z}
    \quad  (\text{resp. }  \langle \lambda, \alpha \rangle \in \mathbb{Z}
    \text{ and } \langle \lambda, \alpha \rangle \geq 0) \quad \text{for all }  \alpha \in \varPhi.
\end{equation}

The set of weights $\Lambda$ forms a subgroup of $\mathcal{E}$ containing the root system $\varPhi$,
i.e. $\varPhi \subset \Lambda \subset \mathcal{E}$.  The concept of a dominant weight was introduced by Cartan
in \cite{C13}, however his definition was differ from the definition \eqref{def_dominant},
see \cite[p.311]{Haw}.

\subsubsection{Partial ordering on the set of weights}
 \label{sec_partial}
Consider two weights $\mu$ and $\lambda$. We say that $\mu$ is {\it higher } than $\lambda$,
and we write $\mu \geq \lambda$ if $\mu - \lambda$ is expressible as a linear  combination
of  positive roots with non-negative real coefficients. This order is only partial.

\begin{proposition}{\cite[ch.VI, $\S$1, $n^{\circ}$6, Prop.18]{Bo02}}
  \label{prop_2}
   The weight $\lambda$ is dominant if and only if
\begin{equation}
     \lambda  \geq w\lambda \; \text{  for any  } \; w \in W.
\end{equation}
\end{proposition}
~\\
The set of dominant vectors is denoted by $\Lambda^+$. We have
\begin{equation}
  \label{eq_domin_weights}
     \Lambda^+ = \Lambda \cap \overline{C}.
\end{equation}

\begin{proposition}
  \label{prop_inique_dominant}
   Any weight is conjugate to unique dominant weight.
\end{proposition}
For details, see \cite[Ch.VI, $\S1$, $n^\circ$10]{Bo02}.
~\\

\subsubsection{Fundamental dominant weights}
  \label{sec_fundam_weights}
The weights $\bar\omega_i$ satisfying  the following relations
\begin{equation}
   \label{fundam_weights}
    \langle \bar\omega_i, \alpha_j \rangle  = \delta_{ij}, \text{ where } i,j \in \{1,\dots,l \}
\end{equation}
are called {\it fundamental dominant weights}. Any weight $\lambda \in \mathcal{E}$ can
be written as an integral linear combination of the vectors $\{\bar\omega_1,\dots,\bar\omega_l \}$.

\subsection{The action of the Weyl group on the the weights}

\subsubsection{Length $l(w)$}
Any element $w$ in the Weyl group $W$ is the product of reflections $s_i$, where
\begin{equation}
     s_i(x) = x - \langle x, \alpha_i \rangle \alpha_i.
\end{equation}

 The minimal number of simple reflections $s_i$ in the decomposition
\begin{equation*}
 w = s_{i_1}\dots s_{i_n}
\end{equation*}
 is called the {\it length} of the element $w$ and is denoted by $l(w)$.
\begin{proposition}{\cite[p.1.7, Corollary]{H90}}
  \label{prop_length_of_w}
  The length of $w$ is equal to the number of positive roots which are transformed
  to negative roots under $w$.
\end{proposition}
The reflection $s_i$ transforms $\alpha_i$ to $-\alpha_i$ and permutes the other
positive roots, then by Proposition \ref{prop_length_of_w}, we have the following

\begin{theorem}{\cite[p.1.6, Lemma]{H90}}
  \label{th_length_plus_1}
 \begin{equation*}
    l(s_iw) =
    \begin{cases}
        & l(w) + 1, \text{ if } w^{-1}(\alpha_i) \in \varPhi^+, \\
        & l(w) - 1, \text{ if } w^{-1}(\alpha_i) \in \varPhi^-.
    \end{cases}
 \end{equation*}
\end{theorem}

\subsubsection{The element of the maximal length}

\begin{proposition}{\cite[ch.VI, $\S$1, $n^o$ 6, Corollary 3]{Bo02}}
 \label{prop_elem_max_len} 
  There exists the unique element $w_0$ of the maximal length in the Weyl group $W$.
  Length $l(w_0)$ is equal to the number of positive roots.
\end{proposition}

\subsubsection{The action $s_i$ on a weight}
  \label{sect_action_si}
Let us expand an arbitrary vector $\lambda \in \mathcal{E}$
in the basis consisting of all fundamental dominant weights $\{\bar\omega_1,\dots, \bar\omega_l \}$:
\begin{equation}
  \label{ex_expansion}
   \lambda = \sum\limits_{i=1}^l m_i\bar\omega_i.
\end{equation}
Here, $(m_1,\dots, m_l)$ are the coordinates of the weight $\lambda$
in the basis $\{\bar\omega_1,\dots, \bar\omega_l\}$.
By \eqref{fundam_weights} we have $m_j = \langle \lambda, \alpha_j \rangle$
for any $j \in \{1,\dots,l \}$.
If $\lambda$ is one of roots, i.e., $\lambda = \alpha_j$ then
\begin{equation}
   \alpha_j = \sum\limits_{i=1}^l c_{ij}\bar\omega_i,
\end{equation}
where $c_{ij} = \langle \alpha_i, \alpha_j \rangle$.
Let $\overline{c}_i = (c_{i1}, \dots, c_{il})$ be the $i$th
row of the Cartan matrix $(\langle \alpha_i, \alpha_j \rangle)_{i,j=1}^l$.
The vector $\overline{c}_i$ is the root $\alpha_i$ in the basis of fundamental weights.
Then
\begin{equation}
  \label{eq_action_si}
  \begin{split}
   & s_i(\lambda) = \lambda - \langle \lambda, \alpha_i \rangle \alpha_i =
        \lambda - m_i(c_{i1},\dots,c_{il}), \text{ i.e., } \\
   & s_i(\lambda) = (m_1 - m_ic_{i1},\dots, m_l - m_ic_{il}).
  \end{split}
\end{equation}
~\\
Eq. \eqref{eq_action_si} is the main formula in Snow's algorithm.

\subsection{Representation and weight space}
  \label{sec_repr_weights}
Let $\mathfrak{g}$ be a Lie algebra over $\mathbb{C}$, and
$\mathfrak{h}$ be a Cartan subalgebra of $\mathfrak{g}$ (a maximal abelian subalgebra).
The roots are defined to be the nonzero eigenvalues of $\mathfrak{h}$ acting on $\mathfrak{g}$
via the adjoint representation:
\begin{equation}
  \alpha : \mathfrak{h} \longrightarrow \mathbb{C},
  \quad [h, x] = \alpha(h)x \text{ for all } h \in \mathfrak{h},
\end{equation}
where $x \in \mathfrak{g}$ is a corresponding eigenvector.
The roots are considered  as linear functionals on $\mathfrak{h}$, they span a real space $E$ in
the dual space $\mathfrak{h}^*$.

Let $V$ be a representation of $\mathfrak{g}$ over $\mathbb{C}$ (not necessarily  adjoint).
A weight $\lambda$ of the representation $V$ with the {\it weight space} of $V_{\lambda}$
is a linear functional on $\mathfrak{h}$  given as follows:
\begin{equation}
  V_{\lambda} := \{x \in V, h\cdot x = \lambda(h) x \text{ for all } h \in \mathfrak{h}\}
\end{equation}

\subsection{Theorem of highest weight}
   \label{sec_highest}
Let $\mathfrak{g}$ be a finite-dimensional semisimple complex Lie algebra.
A weight $\lambda$ of a representation $V$ of $\mathfrak{g}$ is called a {\it highest weight} if
$\mu  \leq \lambda$ for every other weight $\mu$ of $V$, see $\S$\ref{sec_partial}.

In 1913 the theorem of highest weight for representations of simple Lie algebras was completed by E.Caratn.

\begin{theorem}[E.Cartan, \cite{C13}]
  \label{th_Cartan}

{\rm(i)} If $V$ is a finite-dimensional irreducible representation of $\mathfrak{g}$,
then $V$ has a unique highest weight, and this highest weight is dominant integral.

{\rm(ii)} If two finite-dimensional irreducible representations have the same highest weight, they are isomorphic.

{\rm(iii)} For each dominant integral weight $\lambda$ , there exists a finite-dimensional irreducible representation with highest weight $\lambda$.
\end{theorem}

\subsection{\bf Fundamental weights in the case $D_4$}

The dependencies of simple roots $\{ \alpha_1, \alpha_2, \alpha_3, \alpha_4 \}$
and elements of the canonical basis $\{e_1, e_2, e_3, e_4 \}$ are as follows:
\begin{equation}
  \begin{split}
   & \alpha_1 = e_1 - e_2, \quad \alpha_2 = e_2 - e_3,  \quad \alpha_3 = e_3 - e_4, \quad \alpha_4 = e_3 + e_4,
   \\
  & e_1 = \alpha_1 + \alpha_2 + \frac{\alpha_3 + \alpha_4}{2}, \quad
  e_2 = \alpha_2 + \frac{\alpha_3 + \alpha_4}{2}, \\
  & e_3 = \frac{\alpha_3 + \alpha_4}{2}, \quad e_4 = \frac{\alpha_4 - \alpha_4}{2}.
  \end{split}
\end{equation}

By \cite[Table \RomanNumeralCaps{4} ]{Bo02}
the {\it fundamental weights} $\{\bar\omega_1, \bar\omega_2, \bar\omega_3, \bar\omega_4\}$  can be calculated by the following formulas:
\begin{equation}
 \label{calc_fundam}
  \begin{split}
   & \bar\omega_1 = e_1  = \alpha_1 + \alpha_2 + \frac{\alpha_3 + \alpha_4}{2}, \\
   & \bar\omega_2 = e_1 + e_2 = a_1 + 2\alpha_2 + \alpha_3 + \alpha_4, \\
   & \bar\omega_3 = \frac{1}{2}(e_1 + e_2 + e_3 - e_4) = \frac{1}{2}(a_1 + 2\alpha_2 + 2\alpha_3 + \alpha_4), \\
   & \bar\omega_4 = \frac{1}{2}(e_1 + e_2 + e_3 + e_4) = \frac{1}{2}(a_1 + 2\alpha_2 + \alpha_3 + 2\alpha_4). \\
  \end{split}
\end{equation}

Let $B$ denote the Cartan matrix. Then
formulas \eqref{calc_fundam} can also be obtained using the inverse of Cartan matrix $B^{-1}$
as follows:
\begin{equation}
    \bar\omega_i = B^{-1}\alpha_i,
\end{equation}
see \cite[Ch.VI, (14)]{Bo02}.  For the case $D_4$:
\begin{equation} B =
  \begin{bmatrix*}[r]
        2 & -1 & 0 & 0 \\
       -1 & -2 & -1 & -1 \\
        0 & -1 & 2 & 0 \\
        0 & -1 & 0 & 2 \\
 \end{bmatrix*},
 \quad
   B^{-1} =
   \begin{bmatrix*}[r]
        1 &  1 & 1/2 & 1/2 \\
        1 &  2 &  1 &  1 \\
       1/2 & 1 &  1 & 1/2 \\
       1/2 &  1 & 1/2 & 1 \\
 \end{bmatrix*}
\end{equation}

\newpage
\section{\bf The Python implementation}
  \label{app_sect_implem}
\subsection{Root system, generators and number of levels}
  \label{py_root_syst}
The file below ($reflections\_D4.py$) contains information related to the current root system:
reflections, the Cartan matrix and number of levels.
You can change to a different root system, only by modifying this file.
The root system is given as string "D4", "B5", "E6", etc.  The generators of the Weyl group are given as the matrices of the faithful representation. The number of levels is equal to the number of positive roots plus one, see Proposition \ref{prop_length_of_w}.

\begin{python}
'''reflections_D4.py'''
root_system = 'D4'

s1 = [[-1, 1, 0, 0],
      [ 0, 1, 0, 0],
      [ 0, 0, 1, 0],
      [ 0, 0, 0, 1]]
s2 = [[1, 0, 0, 0],
      [1,-1, 1, 1],
      [0, 0, 1, 0],
      [0, 0, 0, 1]]
s3 = [[1, 0, 0, 0],
      [0, 1, 0, 0],
      [0, 1,-1, 0],
      [0, 0, 0, 1]]
s4 = [[1, 0, 0, 0],
      [0, 1, 0, 0],
      [0, 0, 1, 0],
      [0, 1, 0,-1]]
refl = []
refl.append(s1)
refl.append(s2)
refl.append(s3)
refl.append(s4)
''' Cartan matrix '''
Cmatr = \
[[ 2,-1, 0, 0],
 [-1, 2,-1,-1],
 [ 0,-1, 2, 0],
 [ 0,-1, 0, 2]]

''' Number of levels = number of positive roots + 1'''
Nlevels = 13
\end{python}

\subsection{Data structure}
  \label{data_element}
The class {\it Element} contains all the information related to the given element
(and its inverse) of the Weyl group.
Consider, for example, the first element in Table \eqref{tab_level_6_1}
\begin{equation}
 \label{elem_example}
    s_1s_2s_4s_3s_2s_1
\end{equation}
The name of the element \eqref{elem_example} in the class is the following string:
\begin{equation}
    s1.s2.s4.s3.s2.s1
\end{equation}
The information added during
the extended Snow's algorithm is the name of the inverse element, its matrix and its location
in "level". This information is needed to calculate conjugacy classes.

\begin{python}
''' element.py '''
import numpy as np

def isIdentity(M):
    for i in range(len(M)):
        for j in range(len(M[0])):
           if i == j and M[i][j] != 1:
               return False
           elif i!=j and M[i][j] != 0:
               return False
    return True

class Element(object):
    def __init__(self, w, name, name_inv, m, m_inv, n_in_lvl):
        self.weight      = w
        self.name        = name
        self.name_inv    = name_inv
        self.matr        = m
        self.matr_inv    = m_inv
        self.n_in_lvl    = n_in_lvl
        '''we don't know yet the location of inverse element in "level" '''
        self.n_inv_in_lvl = -1

    def keyValAndKeyInv(self):
        key     = ''.join(str(i) for row in self.matr for i in row)
        key_inv = ''.join(str(i) for row in self.matr_inv for i in row)
        val = str(self.n_in_lvl)
        return key, key_inv, val

    def ifSelfInverseMatr(self):
        prod = np.matmul(self.matr, self.matr)
        return isIdentity(prod)
\end{python}

\subsection{Calculation of all levels}
   \label{py_levels}
This section contains the main implementation file of the extended Snow's algorithm
including the search for  inverse elements. To move to another root system,
you need to change the inclusion $from\ reflections\_D4$ to the appropriate one, see Section
\ref{py_root_syst}.

\begin{python}
'''algor_Snow_InvElem.py  '''
import numpy as np
from reflections_D4 import root_system, refl, Cmatr, Nlevels
from element import Element

def buildLevel_0(oneLevel):
   start_weight = np.ones(len_weight, dtype=int).tolist()
   unit_matr = np.eye(len_weight, dtype=int).tolist()
   elm = Element(w=start_weight, name=' ', name_inv = ' ', \
                 m=unit_matr, m_inv=unit_matr, n_in_lvl = 0)
   elm.n_inv_in_lvl = 0
   oneLevel.append(elm)

def newPossibleWeight(numbRefl,  weight, mi):
    new_possible_weight = []
    for jW in range(len_weight):
        ''' The main formula of Snow's algoritm '''
        new_coord =  weight[jW] - mi*Cmatr[jW][numbRefl]
        new_possible_weight.append(new_coord)
    return new_possible_weight

def newElem(iRefl, new_weight, name, name_inv, matr, matr_inv, new_n_in_lvl):
    ''' new_n_in_lvl = the following place in the oneLevel, i.e. = len(one_level) '''
    iW = iRefl - 1
    if (name == ' '):
       new_name_inv = new_name        = str('s') + str(iRefl)
       new_matr_inv = new_matr        = refl[iW]
    else:
       new_name     = str('s') + str(iRefl) + str('.') + name
       new_name_inv = name_inv + str('.s') + str(iRefl)
       new_matr     = np.matmul(refl[iW], matr)
       new_matr_inv = np.matmul(matr_inv, refl[iW])

    new_elem = \
      Element(new_weight, new_name, new_name_inv, new_matr, new_matr_inv, new_n_in_lvl)
    return new_elem

def findAllLevels_to_LvlK(root_system, list_of_all_levels, lvlK):

  len_by_all_levels = 0
  for ik in range(lvlK):
     ''' Get Lvl(k) and create Lvl(k+1) '''
     oneLevel = list_of_all_levels[ik]
     new_level = []
     dictElemsOfLevel =  {}
     len_by_all_levels =  len_by_all_levels + len(oneLevel)

     for iElem in range(len(oneLevel)):
          elem = oneLevel[iElem]
          ''' get elements of lvl = ik to construct the lvl = (ik+1) '''
          weight = elem.weight

          ''' iRefl is  the numb of reflection '''
          for iW in range(len_weight):
            iRefl = iW + 1
            if  weight[iW] > 0:
                mi = weight[iW]
                new_candidate_weight = newPossibleWeight(iW, weight, mi)
                ''' should be unique weight '''
                uniqueFlag = True
                if (iW == len_weight - 1):
                    uniqueFlag = True
                else:
                    for iUniq in range(iW+1, len_weight):
                       if new_candidate_weight[iUniq] < 0:
                           uniqueFlag = False
                           break

                if uniqueFlag is True:
                    ''' This the element of order 2 '''
                    new_n_in_lvl = len(new_level)
                    new_elm = newElem(iRefl, new_candidate_weight, \
                        elem.name, elem.name_inv, \
                        elem.matr, elem.matr_inv, new_n_in_lvl)

                    if new_elm.ifSelfInverseMatr():
                        new_elm.n_inv_in_lvl = new_elm.n_in_lvl
                        ''' no need to save this elem in dictionary'''
                        new_level.append(new_elm)
                    else:
                        key, key_inv, val = new_elm.keyValAndKeyInv()

                        if key in dictElemsOfLevel.keys():
                           ''' the partner (inverse) already waits for this key'''
                           ''' relate 'new_elem_inv' and 'new_elm'  '''
                           val = dictElemsOfLevel[key]
                           n_in_lvl = int(val)
                           new_elem_inv = new_level[n_in_lvl]
                           new_elem_inv.n_inv_in_lvl = new_elm.n_in_lvl
                           new_elm.n_inv_in_lvl = new_elem_inv.n_in_lvl
                           new_level.append(new_elm)
                        else:
                           new_elm.n_in_lvl = len(new_level)
                           ''' inform the partner(inv) about location of new element'''
                           val = str(new_elm.n_in_lvl)
                           dictElemsOfLevel[key_inv] = val
                           new_level.append(new_elm)

     list_of_all_levels.append(new_level)

     ''' write down new_level on a disk '''
     writeOneLevel(ik+1, new_level, prefix)


''' single level recoding procedure. Parameters are as follows:
ik - number of level, oneLevel - one level from list_of_all_levels,
prefix - string root_system, like "D4", "B5", "E6", etc.'''
def writeOneLevel(ik, oneLevel, prefix):

    n_elems = len(oneLevel)

    if n_elems == 0:
        return

    file_name = prefix + '_WeightMatrByLevel_' + str(ik) +\
            '_elems=' + str(n_elems) + '.txt'
    path_to_file = prefix + '_DataFiles\\' + file_name
    print('write file: ', path_to_file)
    with  open(path_to_file, 'w') as f:

        '''in weight the last comma already there '''
        for iElem in range(len(oneLevel)):
            elem = oneLevel[iElem]
            wStr = weightToStr(elem.weight)
            lineElem = 'n='+ str(elem.n_in_lvl) +\
                        ', name=' + elem.name +\
                        ', w=' + wStr +\
                        ', n_inv=' + str(elem.n_inv_in_lvl)
            f.write(lineElem)
            for r in elem.matr:
                line = list(r)
                f.write('\n')
                f.write(str(line))
            f.write('\n')
        f.close()

if __name__ == "__main__":

    list_of_all_levels = []
    len_weight = len(Cmatr)
    oneLevel = []

    ''' Step 0'''
    buildLevel_0(oneLevel)
    writeOneLevel(0, oneLevel, root_system)
    list_of_all_levels.append(oneLevel)

    ''' Here, function writeOneLevel is called for each level  '''
    findAllLevels_to_LvlK(root_system, list_of_all_levels, lvlK=Nlevels)
\end{python}

\subsection{Sample output: file containing one level}
For the root system $D_4$, we get $13$ files corresponding to $13$ levels.
Here is the file containing level $2$ consisting of $9$ elements.
\begin{python}
n=0, name=s2.s1, w=1,-2,3,3, n_inv=3
[-1, 1, 0, 0]
[-1, 0, 1, 1]
[0, 0, 1, 0]
[0, 0, 0, 1]
n=1, name=s3.s1, w=-1,3,-1,1, n_inv=1
[-1, 1, 0, 0]
[0, 1, 0, 0]
[0, 1, -1, 0]
[0, 0, 0, 1]
n=2, name=s4.s1, w=-1,3,1,-1, n_inv=2
[-1, 1, 0, 0]
[0, 1, 0, 0]
[0, 0, 1, 0]
[0, 1, 0, -1]
n=3, name=s1.s2, w=-2,1,2,2, n_inv=0
[0, -1, 1, 1]
[1, -1, 1, 1]
[0, 0, 1, 0]
[0, 0, 0, 1]
n=4, name=s3.s2, w=2,1,-2,2, n_inv=6
[1, 0, 0, 0]
[1, -1, 1, 1]
[1, -1, 0, 1]
[0, 0, 0, 1]
n=5, name=s4.s2, w=2,1,2,-2, n_inv=8
[1, 0, 0, 0]
[1, -1, 1, 1]
[0, 0, 1, 0]
[1, -1, 1, 0]
n=6, name=s2.s3, w=3,-2,1,3, n_inv=4
[1, 0, 0, 0]
[1, 0, -1, 1]
[0, 1, -1, 0]
[0, 0, 0, 1]
n=7, name=s4.s3, w=1,3,-1,-1, n_inv=7
[1, 0, 0, 0]
[0, 1, 0, 0]
[0, 1, -1, 0]
[0, 1, 0, -1]
n=8, name=s2.s4, w=3,-2,3,1, n_inv=5
[1, 0, 0, 0]
[1, 0, 1, -1]
[0, 0, 1, 0]
[0, 1, 0, -1]
\end{python}

\end{appendix}
\newpage
\pagebreak[4]

\end{document}